\RequirePackage{ifpdf}
\documentclass[11pt]{article}
\usepackage{amsthm}
\usepackage{amsmath}
\usepackage{amssymb}
\usepackage{subeqnarray}
\usepackage{cases}
\usepackage{latexsym}
\usepackage{graphicx}
\usepackage{fancybox,color}
\usepackage{longtable}
\usepackage{color}
\usepackage{lscape}
\usepackage{graphicx}
\usepackage{float}
\usepackage[breaklinks=true,bookmarks=true]{hyperref}

\definecolor{lred}{rgb}{1,0.8,0.8}
\definecolor{lblue}{rgb}{0.8,0.8,1}
\definecolor{dred}{rgb}{0.6,0,0}
\definecolor{dblue}{rgb}{0,0,0.5}
\usepackage{url}

\newcommand{\blue}[1]{\begin{color}{blue}#1\end{color}}

\newcommand{\red}[1]{\begin{color}{red}#1\end{color}}
\newcommand{\green}[1]{\begin{color}{green}#1\end{color}}

\def\cA{{\cal A}} 
\def\cB{{\cal B}} \def\tcB{\widetilde{\cB}}
\def\cE{{\cal E}}

\def\cL{{\cal L}}
\def\cD{{\cal D}}
\def\cU{{\cal U}} 

\def\cP{{\cal P}}
 \def\R{{\cal R}}
\def\cS{{\cal S}}
\def\cT{{\cal T}} 
\def\cH{{\cal H}}

\def\cL{{\cal L}}  
\def\cS{{\cal S}} \def\S{{\cal S}} 
\def\cX{{\cal X}}
\def\cY{{\cal Y}}

\def\eps{{\varepsilon}}

\def\grad{\nabla}

\newcommand{\inprod}[2]{\langle #1 , #2 \rangle}

\newcommand{\T}{\mbox{\textrm{\tiny{T}}}}

\def\alp{\alpha}

\def\norm#1{\|#1\|}
\def\abs#1{\left| #1\right|}
\def\cN{{\cal N}} \def\cC{{\cal C}} \def\cD{{\cal D}}

\def\cQ{{\cal Q}} \def\cU{{\cal U}}
 
\def\cK{{\cal K}}
\def\cI{{\cal I}}

\def\hdel{\hat{\delta}}

\def\eps{{\epsilon}}

\def\lam{{\lambda}}

\def\inprod#1#2{\langle#1, \, #2\rangle}

\def\ni{\noindent}

\def\diag{{\rm diag}}

\def\sig{\sigma} 
\def\nn{\nonumber}

 \def\tM{\widetilde{M}}
\def\tS{\widetilde{S}} \def\tZ{\widetilde{Z}} \def\tN{\widetilde{N}}

\def\bP{{\bf P}}
\def\bD{{\bf D}}

\def\domp{\mbox{dom}(p)}
\def\dom{\mbox{dom}}

\def\bF{\bar{F}}

\def\ty{\tilde{y}} 
\def\tv{\tilde{v}} \def\tx{\widetilde{x}}

\newtheorem{thm}{Theorem}[section]  

\newtheorem{prop}[thm]{Proposition}

\newtheorem{rem}{{\it Remark}}           

\DeclareMathOperator*{\argmini}{arg \, min}

\newcommand{\argmin}{\argmini}

\begin{document}

\title{\bf An Efficient Inexact ABCD Method  for  Least Squares Semidefinite Programming }

\author{
	Defeng Sun\thanks{Department  of  Mathematics  and  Risk  Management  Institute, National University of Singapore, 10 Lower Kent Ridge Road, Singapore ({\tt matsundf@nus.edu.sg}). Research supported in part by the Ministry of Education, Singapore, Academic Research Fund under Grant R-146-000-194-112.}, \,
	Kim-Chuan Toh\thanks{Department of Mathematics, National University of Singapore, 10 Lower Kent Ridge Road, Singapore
		({\tt mattohkc@nus.edu.sg}).
		Research supported in part by the Ministry of Education, Singapore, Academic
		Research Fund under Grant R-146-000-194-112.
	}, \,
	Liuqin Yang\thanks{Department of Mathematics, National University of Singapore, 10 Lower Kent Ridge Road, Singapore ({\tt yangliuqin@u.nus.edu}).}
}
\date{May 16,  2015; Revised on May 25,  2015}

\maketitle

\begin{abstract}
We consider least squares semidefinite programming (LSSDP)
where the primal matrix variable must satisfy given linear equality and inequality
constraints, and  must also lie in the intersection of the
cone of symmetric   positive semidefinite matrices and a simple polyhedral set.
We propose an inexact accelerated block coordinate descent (ABCD)
method for solving LSSDP via its dual, which can be reformulated
as a convex composite minimization
problem whose objective is the sum of
a coupled  quadratic  function involving four blocks of variables
and
two separable non-smooth functions involving only the first and second block, respectively. Our inexact  ABCD method has the attractive $O(1/k^2)$ iteration complexity
if the subproblems are solved progressively more accurately.
The design of our ABCD method relies on recent advances
in the symmetric Gauss-Seidel technique for solving a
convex minimization problem whose objective is the sum of a
 multi-block quadratic function  and a
non-smooth function involving only the first block.
Extensive numerical experiments on various classes of over 600 large scale
LSSDP problems demonstrate that our proposed ABCD method not only
can solve the problems to high accuracy, but it is also far more efficient
than (a) the well known BCD (block coordinate descent) method,
(b) the eARBCG (an enhanced version of the accelerated randomized block coordinate gradient) method,
and (c) the APG (accelerated proximal gradient) method.
\end{abstract}

\noindent
\textbf{Keywords:}
Least squares SDP, accelerated block coordinate descent, symmetric Gauss-Seidel.

\medskip
\noindent
\textbf{AMS subject classifications.} 90C06, 90C22, 90C25, 65F10.

\section{Introduction}

Let $\cS^n$ be the space of $n \times n$ real symmetric matrices endowed with the
standard trace inner product $\inprod{\cdot}{\cdot}$ and its induced norm
$\norm{\cdot}$. We denote by  $\S^n_+$ the cone of  positive semidefinite matrices in $\S^n$.  For any matrix $X \in \cS^n$, we use $X \succ 0$
$(X\succeq 0) $ to indicate that $X$ is a symmetric positive definite (positive semidefinite) matrix.

Consider the following semidefinite programming (SDP) problem:
\begin{eqnarray} \label{eq-SDP}
	\begin{array}{ll}	
		\min &  \inprod{C}{X}  \\[5pt]
		{\rm s.t.} &  \cA_E(X) = b_E, \; \cA_I X -s =0, \; X\in \S^n_+, \; X \in \cP, s \in \cK,
	\end{array}
\end{eqnarray}
where $b_E \in \Re^{m_E}$ and $C \in \cS^n$ are given data, $\cA_E:\cX\rightarrow \Re^{m_E}$ and $\cA_I:\cX\rightarrow \Re^{m_I}$ are two given linear maps, $\cP$ and $\cK$ are two nonempty simple closed convex sets,
e.g., $\cP =\{W \in \cS^n:\; L\leq W\leq U\}$ with $L,U\in \cS^n$ being given matrices and $\cK =\{w \in \Re^{m_I}:\; l\leq w\leq u\}$ with $l,u\in \Re^{m_I}$ being given vectors.
When applying a proximal point algorithm (PPA) \cite{rockafellar1976augmented, rockafellar1976monotone} to solve \eqref{eq-SDP}, we need to
solve the following subproblem in each iteration for a given point $(X^k,s^k)$ and a parameter $\sig_k > 0$:
\begin{eqnarray} \label{eq-ppa}
(X^{k+1},s^{k+1}) = \argmin_{X,s}  \left\{
\begin{array}{l}		
\inprod{C}{X} + \frac{1}{2
	\sigma_k}(\norm{X - X^k }^2+\norm{s-s^k}^2)
 \\[5pt]		
\mid \; \cA_E(X) = b_E, \; \cA_I X -s =0, \; X\in \S^n_+, \; X \in \cP, \; s \in \cK
\end{array}\right\}.
\end{eqnarray}
This motivated us to study the
following least squares semidefinite programming (LSSDP)
 which includes \eqref{eq-ppa} as a particular case:
\begin{eqnarray*}
	\begin{array}{lll}	
	\mbox{(\bP)} &	\min &  \frac{1}{2}\norm{X-G}^2 + \frac{1}{2}\norm{s-g}^2  \\[5pt]
		&
{\rm s.t.} &  \cA_E(X) = b_E, \; \cA_I X -s = 0, \; X\in \S^n_+, \; X \in \cP, \; s \in \cK,
	\end{array}
\end{eqnarray*}
where  $G\in \S^n$, $g \in \Re^{m_I}$ are given data.
In order for the PPA to be efficient for solving \eqref{eq-SDP}, it is of great importance for us to design an efficient algorithm to solve the above problem $\mbox{(\bP)}$.
Thus,  the objective of this paper is to achieve this goal via solving the dual of
$\mbox{(\bP)}$.

The dual of $\mbox{(\bP)}$ is given by
\begin{eqnarray*}
	\begin{array}{lrl}	
	\mbox{(\bD)} &	\min F(Z,v,S,y_E,y_I) &:= -\inprod{b_E}{y_E}+ \delta_{\cS^n_+}(S)+ \delta_{\cP}^*(-Z) + \delta_{\cK}^*(-v) \\[5pt]
	&	& + \frac{1}{2} \norm{\cA_E^*y_E+\cA_I^* y _I+S+Z+G}^2 \\[5pt]		
	&	&  +\frac{1}{2}\norm{g+v-y_I}^2
		-\frac{1}{2}\norm{G}^2 -\frac{1}{2}\norm{g}^2,
	\end{array}
\end{eqnarray*}
where for any given  set $\cC$, $\delta_\cC(\cdot)$ is the indicator function over $\cC$ such that $\delta_\cC (u)=0$ if $u\in \cC$ and $\infty$ otherwise,  and $ \delta_{\cC}^*(\cdot) $ is
the conjugate function of $\delta_\cC$  defined by
\begin{eqnarray}
\delta_{\cC}^*(\cdot) =  \sup_{W \in \cC} \inprod{\cdot}{W}. \label{gz}
\end{eqnarray}

Problem  (\bD) belongs to a general class of multi-block convex optimization problems of the form:
\begin{eqnarray} \label{prob-gen}
\min\{ \Psi(x) := \theta(x) + \zeta(x)  \},
\end{eqnarray}
where $x = (x_1,\dots,x_{q}) \in \cX := \cX_1 \times \cdots \times \cX_{q}$,
 and each $ \cX_i$ is a  finite dimensional real Euclidean space equipped with an inner product $\inprod{\cdot}{\cdot}$ and its induced norm $\norm{\cdot}$.
Here
 $\theta(x)=\sum_{i=1}^{q}\theta_i(x_i)$, $\zeta: \cX  \to \R$ and $\theta_i: \cX_i \to (-\infty, +\infty], \, i = 1,\dots, q$ are proper, lower semi-continuous convex functions. We assume that
$\zeta$ is continuously differentiable on an open neighborhood containing $\dom(\theta) := \{ x \in \cX: \theta(x) < \infty\}$   and its gradient
$\grad \zeta$ is  Lipschitz continuous. Note that one  can write ({\bf D}) in the form  of (\ref{prob-gen}) in a number of different ways. One natural choice is of course to
express it in the form of (\ref{prob-gen}) for $q=4$ with $(x_1, x_2, x_3, x_4) \equiv ((Z,v),S, y_E, y_I)$. Another possibility is to express it in the form of
 (\ref{prob-gen}) for $q=2$ with $(x_1, x_2) \equiv ((Z,v), (S, y_E, y_I)).$

For the problem in \eqref{prob-gen},
a well known technique for solving it is the block coordinate descent (BCD) method,  for examples, see \cite{grippo2000convergence, sardy2000block,tseng1993dual,tseng2001convergence} and references therein. Specifically, at iteration $k$, one may update the blocks successively in the Guass-Seidel fashion (other rules can also be applied, see \cite{tseng2001convergence}):
\begin{eqnarray}\label{bcd}
x_1^{k+1} &=& \mbox{$\argmin_{x_1}$}
\Psi(x_1, x_2^k, \ldots, x_q^k), \nonumber
\\[-8pt]
& \vdots & \nonumber
\\[0pt]
x_i^{k+1} &=& \mbox{$\argmin_{x_i}$} \Psi(x_1^{k+1},\ldots, x_{i-1}^{k+1}, x_i, x_{i+1}^k, \ldots, x^k_q),
\\[-8pt]
& \vdots & \nonumber
\\[0pt] \nonumber
x_q^{k+1} &=& \mbox{$\argmin_{x_q}$} \Psi(x_1^{k+1},\ldots,   x_{q-1}^{k+1},  x_q). \nonumber
\end{eqnarray}
When the subproblems in \eqref{bcd} are not easily solvable, a popular approach is to use a single step of the proximal gradient (PG) method, thus yielding the block coordinate gradient descent (BCGD) method \cite{tseng2009coordinate,beck2013convergence}.

Problem \eqref{prob-gen} can also be solved by the accelerated proximal gradient (APG) method with iteration complexity of $O(1/k^2)$ such as in \cite{nesterov1983method,nesterov2004introductory,nesterov2005smooth,beck2009fast,tseng2008accelerated}.
In the best case,  BCD-type methods  have an
iteration complexity of  $O(1/k)$ (see \cite{saha2013nonasymptotic,beck2013convergence}), and can hardly be accelerated to $O(1/k^2)$ as in the case  for the APG method. Nevertheless, some researchers have tried to tackle this difficulty from different aspects. Beck and Tetruashvili \cite{beck2013convergence} proposed an accelerated BCGD method for solving \eqref{prob-gen} by assuming that $\theta \equiv 0$, i.e., without the nonsmooth terms. Very recently, Chambolle and Pock \cite{chambolle2015remark} presented an accelerated BCD method for solving \eqref{prob-gen} by assuming that  $\zeta(x)$ has the special form $\zeta(x) = \sum_{1 \leq i < j \leq q} \norm{A_{ij} x_i + A_{ji} x_j}^2$.
In theory, this method can be applied to the problem (\bD) by choosing
 $x_1 = (Z,v)$ and $x_2 = (S,y_E,y_I)$. But the serious practical disadvantage
is that
the method in \cite{chambolle2015remark}  does not cater for inexact solutions of the associated subproblems and hence
it is not suitable for large scale problems since for (\bD) the subproblems must
be solved numerically.

Besides BCD-type methods based on deterministic updating order, there has been a wide interest in randomized BCD-type methods. Nesterov \cite{nesterov2012efficiency} presented a randomized BCD method with unconstrained and constrained versions in which the selection of the blocks is not done by a deterministic rule (such as the cyclic rule \eqref{bcd}), but rather via a predescribed distribution. Furthermore, an accelerated $O(1/k^2)$ variant was studied for the unconstrained version. To extend the accelerated version for the more generic problem \eqref{prob-gen}, Fercoq and Richt{\'a}rik \cite{fercoq2013accelerated}
proposed an accelerated $O(1/k^2)$ randomized block  coordinate gradient (ARBCG)  method. For strongly convex
functions, Lin, Lu and Xiao \cite{lin2014accelerated} showed that a variant of this
 method can  achieve a linear convergence rate. However, from our numerical experience, the
ARBCG method usually can only solve (\bD) to an accuracy of $10^{-3}$--$10^{-4}$ since only the maximal eigenvalue of $\cA_E \cA_E^*$ is used when  updating $y_E$ (similarly for updating $y_I$). Even a numerically  much  enhanced version of the ARBCG (denoted as eARBCG) method with a weighted norm (for which the theoretical convergence needs to be studied) is also typically $3$--$4$ times slower than the  accelerated  block  coordinate descent  (ABCD) method with a special deterministic rule which we will propose later.

In this paper we aim to design an efficient {\em inexact} accelerated  BCD-type  method
whose worst-case iteration complexity is $O(1/k^2)$
for (\bD). We achieve this  goal by  first proposing an inexact  accelerated  block  coordinate gradient descent  (ABCGD) method  for  the general convex programming problem \eqref{prob-gen} with $\theta_3 = \cdots = \theta_q = 0$, and  then apply it to (\bD) to obtain an
inexact  ABCD  method. Note that
when $\zeta$ is a convex quadratic function,  the ABCGD method and the ABCD method  are identical.
Our ABCD method is designed based on three components. First, we apply a Danskin-type theorem to eliminate the variable $x_1$ in  \eqref{prob-gen}. Then we adapt the
inexact APG framework of Jiang, Sun and Toh proposed in \cite{jiang2012inexact}
to solve the resulting reduced problem. By choosing an appropriate proximal term
and adapting the recently developed inexact symmetric Gauss-Seidel decomposition technique \cite{LiSunToh2015}
for solving a multi-block convex quadratic minimization problem (possibly with a
single nonsmooth block), we show that each subproblem in the inexact APG
method can be solved efficiently in a fashion almost like
the symmetric Gauss-Seidel update.

As already mentioned, one can also apply the APG method directly to solve (\bD), or more generally
\eqref{prob-gen}. In this paper, we also adapt the APG method to directly
solve (\bD) for the sake of numerical comparison. In addition, since the  ARBCG   method does not perform well for solving (\bD), again for the sake of numerical comparison,    we  propose an enhanced version  (called eARBCG) of an accelerated
randomized block coordinate gradient method designed in \cite{lin2014accelerated}
for solving (\bD).  As one  can  see later from the extensive numerical experiments
  conducted to evaluate the performance of various methods,
though the BCD,  APG and eARBCG methods are natural choices for solving (\bD), they are
substantially less efficient than the ABCD method that we have designed.
In particular, for  solving
(\bD),
the ABCD method is  at least  ten times faster than the BCD method for vast majority of the tested problems.
It is quite surprising that a simple novel acceleration step with a special BCD cycle,
as in the case of the ABCD method, can improve the
performance of the standard Gauss-Seidel BCD method by such a
dramatic margin.

The paper is organized as follows. In the next section, we introduce the
key ingredients needed to design our proposed algorithm for solving \eqref{prob-gen},
namely, a Danskin-type theorem for parametric optimization, the inexact
symmetric Gauss-Seidel decomposition technique, and the inexact APG method.
In Section \ref{sec-ABCGD}, we describe the integration of the three ingredients to
design our inexact ABCGD  method for solving \eqref{prob-gen}.
Section \ref{sec-ABCD} presents some specializations of the introduced inexact ABCGD method to
solve the dual LSSDP problem (\bD), as well as discussions on the
numerical computations involved in solving the subproblems.
In Section \ref{sec-APG_eARBCG}, we describe the direct application of the APG method for solving (\bD).
In addition, we also propose an enhancement of the accelerated randomized block
coordinate gradient method in \cite{lin2014accelerated} for solving (\bD).
Extensive numerical experiments to evaluate the performance of the
ABCD, APG, eARBCG and BCD methods are presented in Section \ref{sec-numerical-ABCD}.
Finally, we conclude the paper in the last section.

\bigskip
\noindent
{\bf Notation.}
For any given self-adjoint positive semidefinite operator $\cT$ that maps a real Euclidean space $\cX$ into itself, we use  $ \cT^{1/2}$ to denote  the unique   self-adjoint positive semidefinite operator
such that $\cT^{1/2} \cT^{1/2} =\cT $ and define
\[
\|x\|_{\cT} : = \sqrt{\inprod{x}{\cT x}} =\|\cT^{1/2} x\|\quad \forall\, x\in \cX.
\]

\section{Preliminaries}

\subsection{A Danskin-type theorem}
 Here we shall present a Danskin-type theorem for parametric optimization problems.

 Let $\cX$ and $\cY$ be two  finite dimensional real Euclidean spaces each equipped with an inner product $\inprod{\cdot}{\cdot}$ and its induced norm $\norm{\cdot}$ and $\varphi: \cY \to (-\infty, +\infty] $ be a lower semi-continuous function and $\Omega$ be a nonempty open set in $\cX$.  Denote the effective domain of $\varphi$ by ${\rm dom} (\varphi)$, which is assumed to be nonempty.  Let $\phi(\cdot, \cdot): \cY \times \cX \to (-\infty, +\infty) $ be a  continuous function.
 Define the function $g: \Omega \to [-\infty, +\infty)$  by
 \begin{equation}\label{func:g}
 g(x) = \inf_{y\in {\cal Y}} \{ \varphi(y) +\phi(y,x)\}, \quad x\in \Omega.
 \end{equation}
 For any given $x\in \Omega$, let ${\cal M}(x)$ denote the solution set, possibly empty, to the optimization problem (\ref{func:g}).  The following theorem, largely due to Danskin \cite{Danskin66}, can be proven by essentially  following the proof given in  {\cite[Theorem 10.2.1]{facchinei2003finite}}.

 \begin{thm}\label{thm:Danskin} Suppose that
 	$\varphi: \cY \to (-\infty, +\infty] $ is  a lower semi-continuous function and $\phi(\cdot, \cdot): \cY \times \cX \to (-\infty, +\infty) $ is  a  continuous function.    Assume  that   for every $y\in {\rm dom}(\varphi)\ne \emptyset$,  $ \phi(y,\cdot)$ is  differentiable on $\Omega$ and  $\nabla_x \phi(\cdot, \cdot)$ is continuous on ${\rm dom}(\varphi)  \times \Omega$. Let $x\in \Omega$ be given. Suppose that  there exists an open  neighborhood ${\cal N}\subseteq \Omega$ of $x$ such that for each $x^\prime \in {\cal N}$, ${\cal M}(x^\prime)$ is nonempty and
 	that the set $\cup_{x^\prime\in {\cal N}} {\cal M}(x^\prime)$ is bounded.    Then the following results hold:
 	\begin{description}
 		\item[(i)] The function $g$ is directionally differentiable at $x$ and for any given $d\in \cX$,
 		\[
 		g^\prime(x;d) = \inf_{y\in {\cal M}(x)} \langle \nabla _x \phi(y, x), d\rangle.
 		\]

 		\item[(ii)] If ${\cal M}(x)$ reduces to a singleton, say ${\cal M}(x) = \{y(x)\}$, then $g$ is G$\hat{a}$teaux differentiable at $x$ with
 		\[
 		\nabla g(x) =  \nabla _x \phi(y(x), x) .
 		\]
 	\end{description}
 \end{thm}

 Danskin's Theorem \ref{thm:Danskin} will lead to the following results when convexities on $\varphi$ and $\phi$ are imposed.

 \begin{prop}\label{prop:Danskin-convex}  Suppose that  $\varphi: \cY \to (-\infty, +\infty] $ is a  proper closed  convex function, $\Omega$ is an open   convex set  of ${\cal X}$  and
 	$\phi(\cdot, \cdot): \cY \times \cX \to (-\infty, +\infty) $ is  a  closed  convex function.
 	{}Assume that for every $y\in {\rm dom}(\varphi)$,  $ \phi(y,\cdot)$ is    differentiable  on $\Omega$ and  $\nabla_x \phi(\cdot, \cdot)$ is continuous on ${\rm dom}(\varphi) \times \Omega$ and
 	that for every   $x^\prime \in \Omega$,    ${\cal M}(x^\prime)$ is a singleton, denoted by $\{y(x^\prime)\}$. Then the following results hold:
 	\begin{description}
 		\item[(i)] Let  $x\in \Omega$ be given. If  there exists an   open  neighborhood ${\cal N}\subseteq \Omega$ of $x$  such that    $y(\cdot )$  is  bounded on  any nonempty compact  subset of ${\cal N}$, then the convex function $g$ is   differentiable on  ${\cal N}$ with
 		\[
 		\nabla g(x^\prime) = \nabla_x \phi(y(x^\prime), x^\prime)  \quad \forall \, x^\prime   \in {\cal N}.
 		\]
 		
 		\item[(ii)]     Suppose that there exists a convex    open  neighborhood ${\cal N}\subseteq \Omega$    such that  $y(\cdot )$  is  bounded on  any nonempty compact  subset of ${\cal N}$.
 		Assume that  for any $y\in {\rm dom}(\varphi)$,
 		$\nabla _x\phi(y,\cdot)$ is Lipschitz continuous on ${\cal N}$ and that there
 		exists a self-joint positive semidefinite linear operator  $\Sigma \succeq 0$ such that for all $   x \in {\cal N}$ and $  y\in {\rm dom}(\varphi)$,\[
 		\Sigma\succeq {\cal H} \quad\forall\, {\cal H} \in  \partial  _{xx}^2\phi(y, x) ,
 		\]
 		where for any given $y\in {\rm dom}(\varphi)$, $ \partial_{xx}^2\phi(y,x)$ is the generalized Hessian of $\phi(y,\cdot)$ at $x$.
 		Then
 		\begin{equation}\label{eq:majorization-add}
 		0\le g(x^\prime) -g(x) -\langle \nabla g(x), x^\prime -x\rangle\le \frac{1}{2}\langle  x^\prime -x , \Sigma( x^\prime -x) \rangle \quad \forall\, x^\prime, x\in {\cal N}.
 		\end{equation}
 		Moreover, if  ${\cal N}=\Omega =\cal {X}$, then  $\nabla g(\cdot)$  is Lipschitz continuous on ${\cal X}$ with the Lipschitz constant $||\Sigma \|_2$ (the spectral norm of $\Sigma$)  and  for any $x\in \cX$,
\begin{equation}\label{eq:majorization-add2}
\Sigma \succeq {\cal G} \quad\forall\, {\cal G} \in  \partial  _{xx}^2 g(x),
 \end{equation}
 where  $\partial  _{xx} ^2g(x)$ denotes the generalized Hessian of $g$ at $x$. 
 		\item[(iii)] Assume that  for every $x\in \Omega$, $\phi(\cdot, x)$ is continuously differentiable on an open neighborhood containing ${\rm dom}(\varphi)$. Suppose that
 		 there exist two constants $\alpha>0$ and  $L>0$ such that for all $ x^\prime, x \in \Omega$,
 		\[
 		\langle  \nabla_y \phi(y ^\prime,x)-\nabla_y \phi(y,x), y^\prime-y \rangle\ge \alpha \| y^\prime-y\|^2 \quad
 		\forall\,  y^\prime, y \in {\rm dom}(\varphi),
 		\]
 		and
 		\[
 		\| \nabla_y \phi(y, x^\prime)-\nabla_y\phi(y, x)\|\le L \|x^\prime -x\| \quad
 		\forall\,  y\in {\rm dom}(\varphi).
 		\]
 		Then $y(\cdot)$ is Lipschitz continuous on $\Omega$ such that for all $ x^\prime, x \in \Omega$,
 		\[
 		\|y(x^\prime)-y(x)\|\le (L/\alpha) \|x^\prime -x\|.
 		\]
 	\end{description}
 \end{prop}
 \begin{proof} {\bf Part(i)}. The convexity of $g$ follows from Rockafellar \cite[Theorem 1]{rockafellar1974conjugate} and the rest of part (i) follows from the above Danskin Theorem \ref{thm:Danskin}  and  \cite[Theorem 25.2]{rockafellar1970convex}.
 	
 	{\bf Part (ii)}. For any $x^\prime, x\in  {\cal N}$, we obtain from part (i) and the generalized mean value theorem (MVT) \cite[Theorem 2.3]{hiriart1984generalized}  that
 	\[
 	\begin{array}{lcl}
 	g(x^\prime)-g(x)& =& \displaystyle \phi(y(x^\prime), x^\prime) + \varphi(y(x^\prime)) - [\phi(y(x),x) + \varphi(y(x)) ]\\[2mm]
 	& \le & \displaystyle \phi(y(x), x^\prime) + \varphi(y(x))-[\phi(y(x), x) + \varphi(y(x)) ]
 	\\[2mm]
 	& = & \displaystyle \phi(y(x), x^\prime) - \phi(y(x), x)
 	\\[2mm]
 	& \le &\displaystyle  \langle \nabla_x \phi(y(x), x), x^\prime -x\rangle + \frac{1}{2}\langle  x^\prime -x , \Sigma( x^\prime -x) \rangle \quad \mbox{ (by using the MVT)}
 	\\[2mm]
 	& = &\displaystyle  \langle \nabla g(x), x^\prime -x\rangle +\frac{1}{2}\langle  x^\prime -x , \Sigma( x^\prime -x) \rangle  \quad \mbox{ (by using Part (i))},
 	\end{array}
 	\]
 which, together with the convexity of  $g$, implies that (\ref{eq:majorization-add}) holds.  
 
 Assume that ${\cal N}=\Omega =\cal {X}$.  By  using (\ref{eq:majorization-add}) and  \cite[Theorem 2.1.5]{nesterov2004introductory}, we can assert  that  $\nabla g$ is globally Lipschitz continuous with the Lipschitz constant $\|\Sigma\|_2$. So by Rademacher's Theorem, the Hessian of $g$ exists almost everywhere in $\cX$.  From (\ref{eq:majorization-add}),  we can observe that for any $x\in \cX$ such that the Hessian of  $g$ at $x$ exists, it holds that \[
 \nabla^2_{xx} g(x) \preceq \Sigma. 
 \] 
 	Thus, (\ref{eq:majorization-add2}) follows from the definition of the generalized Hessian of $g$. 
 	
 	{\bf Part (iii)}. The conclusions of part (iii)  follow from part (i),  the maximal monotonicity of $\partial \varphi(\cdot)$, the assumptions  and the fact that for every $x^\prime \in \Omega$,
 	\[
 	0 \in \partial \varphi(y(x^\prime)) +\nabla_y \phi(y(x^\prime), x^\prime).
 	\]
 	We omit the details here.
 \end{proof}

\subsection{An inexact block symmetric  Gauss-Seidel iteration}

Let $s\ge 2$ be a given integer and $\cX := \cX_1\times\cX_2\times\ldots\times\cX_s$,
where $\cX_i$'s are real finite dimensional Euclidean spaces.
 For any $x\in \cX$,   we write $x \equiv (x_1,x_2,\ldots,x_s) $ with $x_i\in \cX_i$.
Let $\cQ : \cX\rightarrow\cX$ be a given self-adjoint positive semidefinite linear operator. Consider the following {block} decomposition
\[\cQ x \equiv \left(
                 \begin{array}{cccc}
                   \cQ_{11} & \cQ_{12} & \cdots & \cQ_{1s} \\
                   \cQ_{12}^* & \cQ_{22} &\cdots  & \cQ_{2s}  \\
                    \vdots& \vdots & \ddots & \vdots \\
                    \cQ_{1s}^*& \cQ_{2s}^* & \cdots & \cQ_{ss} \\
                 \end{array}
               \right)\left(
                        \begin{array}{c}
                          x_1 \\
                          x_2 \\
                          \vdots \\
                          x_s \\
                        \end{array}
                      \right), \quad
\cU x \equiv \left( \begin{array}{cccc}
 0 & \cQ_{12} & \cdots & \cQ_{1s} \\
   & \ddots     &            & \vdots \\
     &                & \ddots & \cQ_{s-1,s} \\
    &                &             &  0
\end{array}\right)\left(\begin{array}{c}
                          x_1 \\
                          x_2 \\
                          \vdots \\
                          x_s \\
                        \end{array}
                      \right),
\]
where $\cQ_{ii}:\cX_i\rightarrow \cX_i, \; i=1,\ldots,s$ are self-adjoint {positive} semidefinite linear operators, $\cQ_{ij}:\cX_j \rightarrow\cX_i, \;i=1,\ldots,s-1 ,\;j>i $ are linear maps.
Note that $\cQ = \cU^* + \cD + \cU$ where $\cD x = (\cQ_{11}x_1,\ldots,\cQ_{ss}x_s)$.

Let $r\equiv (r_1,r_2,\ldots,r_s) \in \cX$ be given.
Define the convex quadratic function $h:\cX\rightarrow\Re$ by
\[h(x): = \frac{1}{2}\inprod{x}{\cQ x} - \inprod{r}{x}, \quad x \in \cX.\]
Let  $p:\cX_1 \rightarrow (-\infty,+\infty]$ be a given lower semi-continuous proper convex function.

Here, we further assume that $\cQ_{ii},\; i=1,\ldots,s$ are positive definite. Define
\[
x_{\leq i} := (x_1,x_2,\ldots,x_i), \quad x_{\geq i} := (x_i,x_{i+1},\ldots,x_s),\; i=0,\ldots,s+1
\] with the convention that
$x_{\leq0} = x_{\geq s+1}=\emptyset.$

Suppose that $\hdel_i, \delta_i^+ \in \cX_i$, $i=1,\ldots, s$ are given error vectors,
with $\hdel_1 = 0$.  Denote $\hdel \equiv (\hdel_1,\ldots,\hdel_s)$ and $\delta^+\equiv(\delta_1^+,\ldots,\delta_s^+)$.
Define the following operator and vector:
\begin{eqnarray}
\cT  &:=& \cU\cD^{-1}\cU^*, \label{eq-O}
\\[5pt]
\Delta(\hdel,\delta^+) &:=& \delta^+ + \cU\cD^{-1}(\delta^+ -\hdel).
\label{eq-Delta}
\end{eqnarray}
Let  $y\in \cX$ be given. Define
\begin{eqnarray}
   x^+ := \argmin_{x}\; \Big\{p(x_1)+h(x)
  +\frac{1}{2}\norm{x - y}_{\cT}^2 -\inprod{\Delta(\hdel,\delta^+)}{x} \Big\}.
 \label{prox-T}
\end{eqnarray}

The following proposition describing an  equivalent BCD-type procedure  for computing $x^+$,
 is the key ingredient for our {subsequent} algorithmic developments.
The proposition
is introduced   by Li, Sun and Toh  \cite{aug2QP}
for the sake of making their  Schur complement   based alternating direction method of multipliers  \cite{LiSunToh2015} more explicit.

\begin{prop}
\label{prop:nbsGS}
 Assume that the self-adjoint linear operators $\cQ_{ii}$,  $i=1,
\ldots,s$ are  positive definite.
Let $y\in \cX$ be given.
For $i=s, \ldots,2,$ define $\hat{x}_i\in \cX_i$ by
\begin{eqnarray}
  \label{y-prime-i}
   \hat{x}_i &:=& \argmin_{x_i} \{\, p(y_1) + h(y_{\le i-1},x_i,\hat{x}_{\ge i+1})-\inprod{\hdel_i}{x_i} \} \nn\\
     &=& \cQ_{ii}^{-1}\big(r_i +\hdel_i - \mbox{$\sum_{j=1}^{i-1}$}
  \cQ_{ji}^*y_j - \mbox{$\sum_{j=i+1}^s$} \cQ_{ij}\hat{x}_j \big).
\end{eqnarray}
 Then the optimal solution $x^+$ defined by \eqref{prox-T} can be obtained exactly {via}
  \begin{eqnarray}\label{prox-nT}
 \left\{ \begin{array}{lcl}
  x_1^+ &=& \argmin_{x_1}\; \{ \, p(x_1) + h(x_1,\hat{x}_{\ge 2}) - \inprod{\delta_1^+}{x_1}\},
\\[5pt]
 x_i^+ &=& \argmin_{x_i} \; \{\, p(x_1^+) + h(x^+_{\le i-1},x_i,\hat{x}_{\ge i+1}) - \inprod{\delta_i^+}{x_i}\}
\\[5pt]
 &=& \cQ_{ii}^{-1}(r_i+\delta_i^+ -\sum_{j=1}^{i-1}\cQ_{ji}^*x_j^+ - \sum_{j=i+1}^s \cQ_{ij}\hat{x}_j),\quad i=2,\ldots,s.
 \end{array} \right.
\end{eqnarray}
Furthermore, $\cH:= \cQ + \cT = (\cD+\cU)\cD^{-1}(\cD+\cU^*)$ is positive definite.
\end{prop}

For later purpose, we also state the following proposition.
\begin{prop} \label{prop-delta}
Suppose that $\cH:= \cQ + \cT = (\cD+\cU)\cD^{-1}(\cD+\cU^*)$ is positive definite.
Let $\xi = \norm{\cH^{-1/2}\Delta(\hdel,\delta^+)} $. Then
\begin{eqnarray}
 \xi = \norm{\cD^{-1/2}(\delta^+-\hdel) + \cD^{1/2}(\cD+\cU)^{-1}\hdel}
\leq \norm{\cD^{-1/2} (\delta^+-\hdel) } +
\norm{\cH^{-1/2} \hdel}.
\end{eqnarray}
\end{prop}
\begin{proof}
The proof is straightforward by using the factorization of $\cH$.
\end{proof}
\begin{rem}
In the computation in \eqref{y-prime-i} and \eqref{prox-nT},
we should interpret the solutions $\hat{x}_i$, $x_i^+$
as  approximate solutions to the minimization problems
without the terms involving $\hdel_i$ and $\delta_i^+$.
Once these approximate solutions have been computed, they would
generate the error vectors $\hdel_i$ and $\delta_i^+$.
With these known error vectors, we know that
the computed approximate solutions are  actually the exact solutions
to the minimization problems in \eqref{y-prime-i} and \eqref{prox-nT}.
\end{rem}

\subsection{An inexact APG method} \label{subsec:APG}
For more generality, we consider the following minimization problem
\begin{eqnarray}
\label{eq-genprob}
\min\{ F(x):= p(x) + f(x)\mid \; x \in \cX \},
\end{eqnarray}
where $\cX$ is a finite-dimensional real Euclidean space.
The functions $f: \cX \rightarrow \Re$,
$p: \cX \rightarrow (-\infty, \infty]$ are proper, lower semi-continuous
convex functions (possibly nonsmooth). We assume that
$f$ is continuously differentiable on $\mathcal{X}$ and its gradient
$\grad f$ is  Lipschitz continuous with modulus $L$ on $\mathcal{X}$, i.e.,
$$
\norm{\grad f(x) - \grad f(y)} \leq L \norm{x-y} \quad \forall\; x,y\in \mathcal{X}.
$$
We also assume that   problem  (\ref{eq-genprob}) is
solvable with an optimal solution $x^*\in \domp$.
The inexact APG algorithm proposed by Jiang, Sun and Toh \cite{jiang2012inexact} for solving \eqref{eq-genprob} is described as follows.

\bigskip
\ni\fbox{\parbox{\textwidth}{\noindent{\bf Algorithm 1.}
		 Input $y^1 = x^0 \in \dom(p)$, $t_1=1 $.
		Set $k=1$. Iterate the following steps.
		\begin{description}
			\item [Step 1.]  Find an approximate minimizer
			\begin{equation}\label{eq-inner-primal}
			x^{k}\approx \mbox{$\argmin_{y\in\cX}$} \Big\{
			f(y^k)+ \inprod{\grad f(y^k)}{y-y^k} + \frac{1}{2}\inprod{y-y^k}{\cH_k(y-y^k)}
			+ p(y) \Big\},
			\end{equation}
			where $\cH_k$ is a self-adjoint positive definite linear operator that is
			chosen by the user.
			\item [Step 2.] Compute $t_{k+1} = \frac{1+\sqrt{1+4t_k^2}}{2}$  and
			$
				y^{k+1} = x^k + \Big(\frac{t_k-1}{t_{k+1}}\Big)(x^k-x^{k-1}).
			$		
		\end{description}}}
		
\bigskip
Given any positive definite linear operator $\cH_j:\cX \rightarrow \cX$, define $q_j(\cdot,\cdot) : \cX \times \cX \rightarrow \Re$ by
\begin{eqnarray*}
q_j(x,y) = f(y) + \inprod{\grad f(y)}{x-y}
+\frac{1}{2} \inprod{x-y}{\cH_j(x-y)}.
\end{eqnarray*}
Let $\{\epsilon_k\}$ be a given convergent sequence
of nonnegative numbers such that
\begin{eqnarray*}
&\sum_{k=1}^{\infty} \epsilon_k < \infty.&
\end{eqnarray*}
Suppose that for each $j$, we have an approximate minimizer:
\begin{eqnarray}
x^j \approx \argmin \{  p(x)+q_j(x,y^j)\mid \; x\in \cX\}
\label{eq-sub}
\end{eqnarray}
that satisfies the following admissible conditions
\begin{eqnarray}
& F(x^j) \;\leq \; p(x^j)+q_j (x^j,y^j),  &
\label{eq-fun-cond}
\\[5pt]
& \grad f(y^j) + \cH_j(x^j -y^j) + \gamma^j = \delta^j \quad
\mbox{with $\norm{\cH_j^{-1/2}\delta^j} \leq \epsilon_j/(\sqrt{2}t_j)$} ,&
\label{eq-grad-cond}
\end{eqnarray}
where $\gamma^j \in \partial p(x^j)$
(Note that for $x^j$ to be an approximate minimizer, we must have
$x^j\in \dom(p)$.)
Then the inexact APG algorithm described in Algorithm 1 has
 the following iteration complexity result.

\begin{thm} \label{thm-apg}
	Suppose that the conditions \eqref{eq-fun-cond} and \eqref{eq-grad-cond} hold and  $\cH_{k-1} \succeq \cH_k \succ 0$ for all $k$.
	Then
	\begin{eqnarray}
	0 \leq F(x^k) - F(x^*) \leq \frac{4}{(k+1)^2}
	\Big( (\sqrt{\tau}+\bar{\epsilon}_k)^2 + 2 \tilde{\epsilon}_k\Big),
	\end{eqnarray}
	where $\tau = \frac{1}{2}\inprod{x^0-x^*}{\cH_1(x^0-x^*)}$, $\bar{\epsilon}_k = \sum_{j=1}^k \epsilon_j$,
	$\tilde{\epsilon}_k = \sum_{j=1}^k \epsilon_j^2.$
\end{thm}
\begin{proof}
See \cite[Theorem 2.1]{jiang2012inexact}.
\end{proof}
\section{An inexact accelerated block coordinate gradient descent method}
\label{sec-ABCGD}

We consider the problem
\begin{eqnarray}
\label{eq-genprob-multiblock}
\min\{ \bF(x_0,x) := \varphi(x_0)+p(x_1) + \phi(x_0,x) \mid \; x_0 \in \cX_0, \; x \in \cX \},
\end{eqnarray}
where $x = (x_1,\dots,x_{s})\in  \cX_1 \times \cdots \times \cX_{s} (=:\cX) $,   $\varphi: \cX_0 \to (-\infty, \infty]$ and $p: \cX_1 \to (-\infty, \infty]$ and
$\phi: \cX_0\times \cX \to \Re$ are three  closed proper convex functions and $\cX_0, \cX_1, \ldots, \cX_s$ are finite dimensional real Euclidean spaces. We assume that
problem  (\ref{eq-genprob-multiblock}) is
solvable with an optimal solution $(x_0^*, x^*)\in \dom(\varphi) \times \dom(p) \times \cX_2 \times \cdots \times \cX_{s}$. Define the function $f: \cX  \to [-\infty, +\infty)$  by
\begin{equation}\label{func:barf}
f(x) = \inf_{x_{0} \in \cX_{0}} \{ \varphi(x_{0}) + \phi(x_0,x)\}, \quad x = (x_1,\dots,x_s) \in \cX.
\end{equation}
For any given $x\in \cX$, let ${\cal M}(x)$ denote the solution set to the optimization problem in \eqref{func:barf}. Suppose that  the assumptions  in Part (ii)\footnote{For subsequent discussions, we do not need the assumptions in  Part (iii) of  Proposition \ref{prop:Danskin-convex} though the results there  have their own merits.}  of Proposition \ref{prop:Danskin-convex} imposed on $\varphi$ and $\phi$  hold for $\Omega =\cX$.
Then, by  Part (ii) of Proposition \ref{prop:Danskin-convex}, we know that   $f$ is continuously differentiable and its gradient $\grad f$ is  globally Lipschitz continuous, and for all $x, y\in \cX$,
\begin{equation}\label{eq:hahaha}
f(x) -f(y) -\langle \nabla f(y), x -y\rangle\le\frac{1}{2}\langle  x -y , \cL( x -y) \rangle  ,
\end{equation}
where $\cL:\cX \rightarrow \cX$ is a self-adjoint positive semidefinite linear operator
such that for any $ x \in \cX$,
\begin{eqnarray}
\cL \succeq \cH  \quad \forall\,  \cH \in \partial^2 f(x),
\end{eqnarray}
where $\partial^2 f(x)$ is the generalized  Hessian   of $f$ at $x$. One natural choice, though not necessarily the best,  for
$\cL $ is $\cL =\Sigma$, where  $\Sigma$ is a self-adjoint positive semidefinite linear operator that satisfies
\begin{eqnarray}
\Sigma \succeq \cH  \quad \forall\,  \cH \in \partial^2_{xx} \phi(x_0,x) \quad \forall \; x_0 \in {\rm dom} \varphi, \; x \in \cX.
\end{eqnarray}
Here for any given $x_0 \in {\rm dom} \varphi$, $\partial^2_{xx} \phi(x_0,x)$ is the generalized  Hessian   of $\phi(x_0,\cdot)$ at $x$.

Now we consider an equivalent problem of \eqref{eq-genprob-multiblock}:
\begin{eqnarray}
\label{eq-genprob1}
\min\{ F(x):= p(x_1) + f(x) \mid \; x \in \cX \}.
\end{eqnarray}
Given a self-adjoint positive semidefinite linear operator $\cQ:\cX \rightarrow \cX$ such that
\begin{eqnarray*}
	\cQ \succeq \cL, \; \cQ_{ii} \succ 0, \; i = 1,\dots,s.
\end{eqnarray*}
Define $\cT$ and $\Delta$ by \eqref{eq-O} and \eqref{eq-Delta}, respectively, and
$h(\cdot,\cdot) : \cX \times \cX \rightarrow \Re$ by
\begin{eqnarray}
h(x,y) &:=& f(y) + \inprod{\grad f(y)}{x-y}
+\frac{1}{2} \inprod{x-y}{\cQ(x-y)}, \nn
\\
&=& f(y) + \inprod{\grad_x \phi(x_0(y),y)}{x-y}
+\frac{1}{2} \inprod{x-y}{\cQ(x-y)}, \quad x, \; y \in \cX,
\end{eqnarray}
where $x_0(y)$  is the unique solution to $\inf_{x_{0} \in \cX_{0}} \{ \varphi(x_{0}) + \phi(x_0,y)\}$. From (\ref{eq:hahaha}),
we have
\begin{eqnarray} \label{eq-majorization}
f(x) \leq h(x,y) \quad \forall x, \; y \in \cX.
\end{eqnarray}

We can now apply the inexact APG method described in  Algorithm 1  to problem \eqref{eq-genprob1} to obtain the following inexact accelerated block coordinate gradient descent (ABCGD) algorithm for problem \eqref{eq-genprob-multiblock}.

\bigskip
\ni\fbox{\parbox{\textwidth}{\noindent{\bf Algorithm 2.}
		 Input $y^1 = x^0 \in \dom(p) \times \cX_2 \times \cdots \times \cX_s$, $t_1=1 $. Let $\{\epsilon_k\}$ be a  summable sequence of nonnegative numbers.
		Set $k=1$. Iterate the following steps.
		\begin{description}
			\item [Step 1.] Suppose
            $\delta_i^k, \, \hat{\delta}_i^k \in \cX_i$, $i=1,\ldots, s$, with $\hat{\delta}_1=0$,
are  error  vectors
such that
			\begin{eqnarray} \label{eq-deltacond}
			  \max\{ \norm{\delta^k}, \norm{\hat{\delta}^k}\} \leq \epsilon_k/(\sqrt{2}t_k).
			\end{eqnarray}			
			Compute
			\begin{eqnarray}
			x^{k}_{0} &=& \argmin_{x_{0}} \Big\{ \varphi(x_{0}) + \phi(x_0,y^k)\Big\},
			\\	
			x^k &=& \argmin_{x} \Big\{p(x_1)+h(x,y^k)
			+\frac{1}{2}\norm{x - y^k }_{\cT}^2 - \inprod{\Delta(\hat{\delta}^k,
\delta^k)}{x} \Big\}.	
			\label{eq-x-update}				
			\end{eqnarray}
			\item [Step 2.] Compute $t_{k+1} = \frac{1+\sqrt{1+4t_k^2}}{2}$ and
			$
				y^{k+1} = x^k + \Big(\frac{t_k-1}{t_{k+1}}\Big)(x^k-x^{k-1}).
			$
		\end{description}}}
		
\bigskip

 The iteration complexity result for the inexact ABCGD algorithm described in
Algorithm 2 follows from Theorem \ref{thm-apg} without much difficulty.
\begin{thm} \label{thm-abcgd}
	Denote $\cH = \cQ+\cT$ and $\beta = 2\norm{\cD^{-1/2}} +\norm{\cH^{-1/2}}$.
	Let $\{(x_0^k,x^k)\}$ be the sequence
	generated by Algorithm 2.
	Then
	\begin{eqnarray} \label{eq-abcgd-result}
	0 \leq F(x^k) - F(x^*) \leq \frac{4}{(k+1)^2}
	\Big( (\sqrt{\tau}+\beta\bar{\epsilon}_k)^2 + 2 \beta^2\tilde{\epsilon}_k\Big),
	\end{eqnarray}
	where $\tau = \frac{1}{2}\inprod{x^0-x^*}{\cH(x^0-x^*)}$, $\bar{\epsilon}_k = \sum_{j=1}^k \epsilon_j$,
	$\tilde{\epsilon}_k = \sum_{j=1}^k \epsilon_j^2.$ Consequently,
	\begin{eqnarray}
	0 \leq \bF(x_0^k,x^k) - \bF(x_0^*,x^*) \leq \frac{4}{(k+1)^2}
	\Big( (\sqrt{\tau}+\beta\bar{\epsilon}_k)^2 + 2 \beta^2\tilde{\epsilon}_k\Big).
	\end{eqnarray}
\end{thm}

\begin{proof}
	From \eqref{eq-majorization}, we have
	\begin{eqnarray}
	F(x^k) \leq p(x_1^k)+h(x^k,y^k)
	+\frac{1}{2}\norm{x^k - y^k }_{\cT}^2.
	\end{eqnarray}
	From the optimality condition for \eqref{eq-x-update}, we can obtain
	\begin{eqnarray}
	\grad_x \phi(x_0^k,y^k) + \cH(x^k -y^k) + \gamma^k = \Delta( \hdel^k, \delta^k ),
	\end{eqnarray}
	where $\gamma^k \in \partial p(x_1^k )$.  By using \eqref{eq-deltacond} and Proposition \ref{prop-delta}, we know that
	\begin{eqnarray*}
    &&\norm{\cH^{-1/2}\Delta(\hdel^k,\delta^k)} \leq
   \norm{\cD^{-1/2}(\delta^k-\hdel^k)} + \norm{\cH^{-1/2}\hdel^k}
    \\[5pt]
    &&\leq\; (2\norm{\cD^{-1/2}} + \norm{\cH^{-1/2}}) \max \{\norm{\delta^k},\norm{\hdel^k}\}
    \;\leq\;  \frac{\beta \epsilon_k}{\sqrt{2}t_k}.
	\end{eqnarray*}
	Then \eqref{eq-abcgd-result} follows from Theorem \ref{thm-apg}.
\end{proof}
				
\begin{rem}
(a) Note that we can use the symmetric  Gauss-Seidel iteration described in Proposition \ref{prop:nbsGS} to compute $x^k$ in \eqref{eq-x-update}. Therefore, Step 1 is actually one special   block coordinate gradient descent (BCGD) cycle with the order $x_0 \rightarrow \hat{x}_s \rightarrow \hat{x}_{s-1} \rightarrow \cdots \rightarrow \hat{x}_2
 \rightarrow x_1 \rightarrow x_2 \cdots \rightarrow x_s$.
\\[5pt]
(b) Assuming that
$x_0^k,\hat{x}_s^k,\dots,\hat{x}_2^k$ and $x_1^k$ have been computed.
One may try to estimate $x_2^k,\dots,x_s^k$ by using $\hat{x}^k_2,\dots,\hat{x}^k_s$,
respectively. In this case the corresponding residual vector $\delta^k$ would be given by:
\begin{eqnarray}
 & \delta^k_i = \hat{\delta}^k_i + \sum_{j=1}^{i-1} \cQ_{ji}^* (\hat{x}^k_i-y^k_i),\quad
i=2,\dots,s. &
\end{eqnarray}
In practice,  we may accept such an approximate solution $x^k_i = \hat{x}^k_i$ for $i=2,\dots,s$, provided  a condition of the form $\norm{\delta^k_i} \leq c \norm{\hat{\delta}^k_i}$  (for some
constant $c > 1$ say $c=10$) is satisfied for $i=2,\dots,s$.
When such an approximation is admissible, then
the linear systems involving $\cQ_{ii}$ need only be solved once
instead of twice for $i=2,\dots,s$. We should emphasize that such a saving can be
significant when the linear systems are solved by a Krylov iterative method
such as the preconditioned conjugate gradient method. Of course, if the linear systems  are
solved by computing the Cholesky factorization (which is done
only once at the start of
the algorithm) of $\cQ_{ii}$, then the
saving would not be as substantial.
\end{rem}
				
\section{Two variants of the inexact ABCD  method for (D)}
\label{sec-ABCD}

Now we can apply Algorithm 2 directly to (\bD) to obtain two variants of the inexact accelerated block coordinate descent method. In the first variant,  we apply Algorithm 2 to (\bD)
by expressing it in the form of \eqref{prob-gen}
with $q=4$ and $(x_1, x_2, x_3, x_4) = ((Z,v),S, y_E, y_I)$.
In the second variant, we express (\bD)
 in the form of
 \eqref{prob-gen} with $q=2$ and $(x_1, x_2) = ((Z,v), (S, y_E, y_I)).$  For the remaining parts of this paper, we  assume that $\cA_E: \cX \to \Re^{m_E}$ is onto and the solution set of (\bD) is nonempty.
 The convergence  of both variants follows from  Theorem  \ref{thm-abcgd} for Algorithm 2 and will not be repeated here.

\bigskip
The detailed steps of the first variant of the inexact ABCD method are given as follows.
\centerline{\fbox{\parbox{\textwidth}{
			{\bf Algorithm ABCD-1}: {\bf An inexact ABCD method for (\bD).}
			\\[5pt]
			Select an initial point  $(\tZ^1,\tv^1,\tS^1,\ty_E^1,\ty_I^1) = (Z^{0},v^0,S^{0},y_E^{0},y_I^0)$ with $(-Z^{0},-v^0)\in
\dom(\delta^*_\cP) \times \dom(\delta^*_\cK)$. Let $\{\epsilon_k\}$ be a  summable sequence of
nonnegative numbers, and $t_1=1 $.
			Set $k=1$. Iterate the following steps.
			\begin{description}				
				\item[Step 1.]  Suppose that $\delta_E^k, \, \hat{\delta}_E^k \in \R^{m_E}$ and  $\delta_I^k, \, \hat{\delta}_I^k \in \R^{m_I}$  are error vectors such that
				\begin{eqnarray*}
				\max \{ \norm{\delta_E^k}, \norm{\delta_I^k},  \norm{\hat{\delta}_E^k}, \norm{\hat{\delta}_I^k} \}\leq \epsilon_k/(\sqrt{2}t_k).
				\end{eqnarray*}
				Let $\widetilde{R}^k = \cA_E^*\ty_E^{k}+\cA_I^*\ty_I^{k}+\tS^{k}+G$. Compute
				\begin{eqnarray*}
			\hspace{-0.7cm}
   	(Z^{k},v^k) &=& \mbox{$\argmini_{Z,v}$}
           \big\{F(Z,v,\tS^{k},\ty_E^{k},\ty_I^{k})\big\}				
				=(\Pi_{\cP}(\widetilde{R}^k)-\widetilde{R}^k, \; \Pi_{\cK}(g-\ty_I^k)-(g-\ty_I^k)),			
				\\[5pt]
				\hat{y}_E^{k} &=& \mbox{$\argmini_{y_E}$}
     \big\{F(Z^k,v^k,\tS^k,y_E,\ty_I^{k})-\inprod{\hat{\delta}_E^k}{y_E}\big\},	
				\\[5pt]
				\hat{y}_I^{k} &=& \mbox{$\argmini_{y_I}$}
     \big\{F(Z^k,v^k,\tS^k,\hat{y}_E^{k},y_I)-\inprod{\hat{\delta}_I^k}{y_I}\big\},				
				\\[5pt]
				S^{k} &=& \mbox{$\argmini_{S}$}
 \big\{F(Z^k,v^k,S,\hat{y}_E^{k},\hat{y}_I^{k})\big\}				
				\;=\;
\Pi_{\S^n_+}\Big(-(\cA_E^*\hat{y}_E^{k}+\cA_I^*\hat{y}_I^{k}+Z^{k}+G)\Big),
				\\[5pt]
				y_I^{k} &=& \mbox{$\argmini_{y_I}$}
 \big\{F(Z^k,v^k,S^k,\hat{y}_E^{k},y_I)-\inprod{\delta_I^k}{y_I}\big\},				
				\\[5pt]
				y_E^{k} &=& \mbox{$\argmini_{y_E}$}
 \big\{F(Z^k,v^k,S^k,y_E,y_I^{k}) - \inprod{\delta_E^k}{y_E}\big\}.
				\end{eqnarray*}
				\item [Step 2.] Set $t_{k+1} = \frac{1+\sqrt{1+4t_k^2}}{2}$ and
    $\beta_k = \frac{t_k-1}{t_{k+1}}$. Compute
				\begin{eqnarray*}
				\hspace{-0.7cm}
      	\tS^{k+1} = S^k + \beta_k(S^k-S^{k-1}),\;					
					\ty_E^{k+1} = y_E^k + \beta_k(y_E^k-y_E^{k-1}), \;				
					\ty_I^{k+1} = y_I^k + \beta_k(y_I^k-y_I^{k-1}).				
				\end{eqnarray*}							
				\end{description}
				}}}
				\bigskip
				
\begin{rem}
(a)
To compute $\hat{y}_E^k$ in Algorithm ABCD-1, we need to solve the following linear system of equations:
	\begin{eqnarray} \label{eq-linear-equality}	
	(\cA_E \cA_E^*)\hat{y}_E^{k} \;\approx\;
 b_E-\cA_E(\cA_I^*\ty_I^{k}+\tS^k+Z^k+G)
	\end{eqnarray}
	with the residual norm
	\begin{eqnarray} \label{eq-linear-equality-delta}	
	\norm{\hat{\delta}_E^k} = \norm{ b_E-\cA_E(\cA_I^*\ty_I^{k}+\tS^k+Z^k+G)
- \cA_E\cA_E^* \hat{y}^k} \leq \epsilon_k / (\sqrt{2} t_k).
	\end{eqnarray}
If the (sparse) Cholesky factorization of $\cA_E \cA_E^*$  can be computed (only once) at a moderate cost, then \eqref{eq-linear-equality} can be solved exactly (i.e., $\hat{\delta}_E^k = 0$); otherwise \eqref{eq-linear-equality} can be solved inexactly to satisfy \eqref{eq-linear-equality-delta} by an iterative method such as the preconditioned conjugate gradient (PCG) method. The same remark also applies to the computation of $\hat{y}_I^{k}, y_I^k, y_E^k$.
\\[5pt]
(b)
From the presentation in Step 1 of Algorithm ABCD-1, it appears that we would need to solve
the linear systems involving the matrices $\cA_E\cA_E^*$ and $\cA_I\cA_I^* + \cI$ twice.
In practice, one can often avoid solving the linear systems twice if
$\hat{y}_E^k$ and $\hat{y}_I^k$ are already sufficiently accurate approximate solutions for the
respective second linear systems. More specifically, suppose that we approximate $y_I^k$ by $\hat{y}_I^k$. Then the residual vector for the second linear system would be given by
$$
 \delta_I^k = \hdel_I^k + \cA_I(S^k - \tS^k).
$$
If the condition
$
 \norm{\delta_I^k} \leq \eps_k/(\sqrt{2}\,t_k)
$
is satisfied, then we need not solve the second linear system since
$\hat{y}_I$ is already a sufficiently accurate solution for the
second linear system. Similarly, if we use $\hat{y}_E^k$ to approximate $y_E$,
then the corresponding residual vector would be given by
$$
 \delta_E^k = \hdel_E^k + \cA_E (\cA_I^* (y_I^k-\ty_I^k)+S^k-\tS^k).
$$
Again if the condition that $\norm{\delta_E^k} \leq \eps_k/(\sqrt{2}t_k)$ is satisfied,
then we can take $y_E^k = \hat{y}_E^k$.
\end{rem}

For  the second variant of the inexact ABCD method, we apply
Algorithm 2 to (\bD) by expressing it
 in the form of
 \eqref{prob-gen} with $q=2$ and $(x_1, x_2) = ((Z,v), (S, y_E, y_I)).$
In this case, we treat $(S,y_E,y_I)$ as a single block and the corresponding subproblem
in the ABCD method neither admits an analytical solution nor is solvable via a linear system
of equations. To solve the  subproblem, we use a semismooth Newton-CG (SNCG) algorithm introduced in \cite{zhao2010newton, yang2014sdpnal} to solve it inexactly.

\bigskip
The detailed steps of the second variant of the inexact ABCD method are given as follows.
We should mention that it is {in} fact an accelerated version of a majorized semismooth Newton-CG  (MSNCG) algorithm presented in \cite{yang2014sdpnal}.

\bigskip
\centerline{\fbox{\parbox{\textwidth}{
			{\bf Algorithm ABCD-2}: {\bf An inexact ABCD method with SNCG for (D).}
			\\[5pt]
			Select an initial point  $(\tZ^1,\tv^1,\tS^1,\ty_E^1,\ty_I^1) = (Z^{0},v^0,S^{0},y_E^{0},y_I^0)$  with $(-Z^{0},-v^0)\in
\dom(\delta^*_\cP) \times \dom(\delta^*_\cK)$. Let $\{\epsilon_k\}$ be a nonnegative summable sequence, $t_1=1 $ and $\tau =10^{-6}$.
			Set $k=1$. Iterate the following steps.
			\begin{description}				
				\item[Step 1.]  Suppose   $\delta_E^k \in \R^{m_E}$,  $\delta_I^k \in \R^{m_I}$  are error vectors such that
				$$\max \{ \norm{\delta_E^k},  \norm{\delta_I^k}\}\leq \epsilon_k/(\sqrt{2}t_k).$$
Let $\widetilde{R}^k = \cA_E^*\ty_E^{k}+\cA_I^*\ty_I^{k}+\tS^{k}+G$.
				Compute
				\begin{eqnarray*}
				&&	\hspace{-0.7cm}
 (Z^{k},v^k) = \argmini_{Z,v} \big\{F(Z,v,\tS^{k},\ty_E^{k},\ty_I^{k})\big\} =(\Pi_{\cP}(\widetilde{R}^k)-\widetilde{R}^k, \Pi_{\cK}(g-\ty_I^k)-(g-\ty_I^k)),
					\\
			&&	\hspace{-0.7cm}
(S^k,y_E^{k},y_I^k) = \argmini_{S,y_E,y_I}
\left\{
F(Z^k,v^k,S,y_E,y_I) +\frac{\tau}{2}\|y_E-\ty_E^k\|^2
-\inprod{\delta_E^k}{y_E}-\inprod{\delta_I^k}{y_I}
\right\}.					
				\end{eqnarray*}				
				\item [Step 2.] Set $t_{k+1} = \frac{1+\sqrt{1+4t_k^2}}{2}$, $\beta_k=\frac{t_k-1}{t_{k+1}}$.
               Compute
				\begin{eqnarray*}
				\hspace{-0.7cm}
               	\tS^{k+1} = S^k +\beta_k (S^k-S^{k-1}), \;					
					\ty_E^{k+1} = y_E^k + \beta_k(y_E^k-y_E^{k-1}),\;					
					\ty_I^{k+1} = y_I^k + \beta_k (y_I^k-y_I^{k-1}).	
				\end{eqnarray*}					
			\end{description}
		}}}

\bigskip
In our numerical experiments, we always start with the first variant of the
ABCD method, and then switch it to the second variant when the convergence speed
of the first variant is deemed to be unsatisfactory.
 {As discussed in \cite{zhao2010newton,yang2014sdpnal}, each iteration of the SNCG algorithm can be quite expensive. In fact, ABCD-1 can achieve a high accuracy efficiently for  most of the  problems perhaps because  it has $O(1/k^2)$ iteration complexity and    there is no need to be switched to ABCD-2. However, for some difficult problems, ABCD-1 may stagnate. In this case, ABCD-2 can perform much better since it wisely makes use of second-order information and it has less blocks.}

\subsection{An efficient iterative  method for solving the linear systems}

Observe that in both the ABCD method and the APG method to be presented later in the
next section, for solving (\bD), we need to solve the following linear systems
\begin{eqnarray} \label{eq-linearsystem}
\cB y_I = r,
\end{eqnarray}
where $\cB = (\cA_I \cA_I^* + \alpha \cI )$, with $\alpha = 1$ and $\alpha = \frac{1}{3}$ for the ABCD method and APG method, respectively.
For the case where the matrix $\cB$ and
its (sparse) Cholesky factorization  (need only to be done once) can be computed at a moderate cost,
\eqref{eq-linearsystem} can be solved efficiently.
However, if the Cholesky factorization of $\cB$ is not available, then we need an alternative efficient method to deal with \eqref{eq-linearsystem}.
In this paper,  we solve \eqref{eq-linearsystem} by a preconditioned conjugate gradient (PCG) method. In order to speed up the convergence of the CG method, we construct the following preconditioner $\widetilde{\cB}$ based {on} a few leading eigen-pairs of $\cB$. Specifically,
consider the following eigenvalue decomposition:
\begin{eqnarray}\label{decomp-AAt}
\cB = PDP^T,
\end{eqnarray}
where $P \in\Re^{m_I\times m_I} $ is an orthogonal matrix whose columns are eigenvectors of $\cB$,  and $D$ is the corresponding diagonal matrix of eigenvalues
with the diagonal elements arranged in a descending order:
$\lambda_1 \geq \lambda_2 \geq \cdots \geq \lambda_{m_I}$. We choose
the preconditioner to be
\begin{eqnarray}
\widetilde{\cB} = \sum_{i=1}^{k_I}\lambda_i P_i P_i^T +\lambda_{k_I}\sum_{i=k_I+1}^{m_I} P_i P_i^T,
\label{eq-tcB}
\end{eqnarray}
where $P_i$ is the $i$-th column of $P$,  $k_I$ is a small integer such that  $1 \leq k_I \ll m_I$ and $\lambda_{k_I} > 0$. Note that we only need to compute
 (which only needs to be done once) the first $k_I$ eigenvalues of $\cB$ and their corresponding eigenvectors. Then
\begin{eqnarray*}
	\tcB^{-1} &=& \sum_{i=1}^{k_I}\lambda_i^{-1} P_i P_i^T +\lambda_{k_I}^{-1}\sum_{i=k_I+1}^{m_I} P_iP_i^T
	\;=\; \sum_{i=1}^{k_I}\lambda_i^{-1} P_i P_i^T +\lambda_{k_I}^{-1}(\cI-\sum_{i=1}^{k_I} P_i P_i^T)
	\\
	&=&\lambda_{k_I}^{-1}\cI-\sum_{i=1}^{k_I-1}(\lambda_{k_I}^{-1} - \lambda_i^{-1}) P_i P_i^T.
\end{eqnarray*}
From the expression of $\tcB$, we can see that the overhead cost of applying the
preconditioning step $\tcB^{-1} v$ for a given $v$ can be kept low
compared to the cost of performing $\cB v$.

In our numerical experiments, we solve \eqref{eq-linearsystem} approximately
by applying the PCG method with the preconditioner \eqref{eq-tcB}
whenever the sparse Cholesky factorization of $\cB$ is too expensive to be computed.
As we are solving a sequence of linear systems of the form \eqref{eq-linearsystem}
where the right-hand side vector changes moderately from iteration to iteration,
we can warm-start the PCG method with the previous solution $y_I^{k-1}$ when
solving the $k$th linear system.
For the problems which we have tested in our numerical experiments,
where the number of linear inequality constraints $m_I$ in (\bP) can be very large,
usually we only need less than ten PCG iterations on the average
to solve \eqref{eq-linearsystem}  to  the required accuracy.
This confirms that the preconditioner \eqref{eq-tcB} is quite effective
for the problems under our consideration, and we have thus presented an efficient iterative method to solve the large scale linear systems here.

\section{An APG method and an enhanced ARBCG method for solving (D)}
\label{sec-APG_eARBCG}

\def\hF{\widehat{F}}

Instead of the ABCD method, one can also apply the APG
and accelerated randomized block coordinate gradient descent
methods to solve (\bD).
The details are given in the next two subsections.

\subsection{An APG method for solving (D)}

To apply the APG method, we note that  (\bD) can equivalently be rewritten as follows:
\begin{eqnarray*}	
	\begin{array}{lll}	
	\mbox{($\widehat{\bD}$)} 	\quad
   \min\;  \hF(y_E,y_I,Z)
:= & -\inprod{b_E}{y_E}+ \delta_{\cP}^*(-Z)   + \frac{1}{2} \norm{\Pi_{\S^n_+}(\cA_E^*y_E+\cA_I^* y _I+Z+G)}^2	
     \\[5pt]			
		&  +\frac{1}{2}\norm{g-y_I}^2 - \frac{1}{2}\norm{(g-y_I)-\Pi_{\cK}(g-y_I)}^2
		-\frac{1}{2}\norm{G}^2 -\frac{1}{2}\norm{g}^2.
	\end{array}
\end{eqnarray*}
In order to apply the APG method to solve \mbox{($\widehat{\bD}$)},
we first derive a majorization of the objective function given
the auxiliary iterate $(\ty_E^k,\ty_I^k,\tZ^k)$. Let $\varphi_1(M,N,Z) = \frac{1}{2}\norm{\Pi_{\cS^n_+}(M+N+Z+G)}^2$,  $\varphi_2(y_I) = \frac{1}{2}\norm{g-y_I}^2 - \frac{1}{2}\norm{(g-y_I)-\Pi_{\cK}(g-y_I)}^2$. It is known that $\grad \varphi_i(\cdot), \, i =1,2$ are Lipschitz continuous.
Thus we have that
\begin{eqnarray*}
	\varphi_1(M,N,Z) &\leq & \varphi_1(\tM^k,\tN^k,\tZ^k) + \inprod{\Pi_{\S^n_+}(\tM^k+\tN^k+\tZ^k+G)}{M+N+Z-\tM^k-\tN^k-\tZ^k} \\
	&& + \frac{3}{2}\norm{M-\tM^k}^2 + \frac{3}{2}\norm{N-\tN^k}^2 + \frac{3}{2} \norm{Z-\tZ^k}^2,
	\\[5pt]
	\varphi_2(y_I) &\leq & \varphi_2(\ty_I^k) - \inprod{\Pi_{\cK}(g-\ty_I^k)}{y_I-\ty_I^k} + \frac{1}{2}\norm{y_I-\ty_I^k}^2,
\end{eqnarray*}
where $M = \cA_E^*y_E$, $N =  \cA_I^*y_I$, $\tM^k = \cA_E^* \ty_E^k$,
$\tN^k = \cA_I^* \ty_I^k$.
From here, we get
\begin{eqnarray*}
	& & \hF(y_E,y_I,Z) - \hF(\ty_E^k,\ty_I^k,\tZ^k) \\
	& \leq& \delta_{\cP}^*(-Z) - \inprod{b_E}{y_E-\ty_E^k} +\inprod{X^k}{\cA_E^*y_E-\cA_E^*\ty_E^k} +\inprod{X^k}{\cA_I^*y_I-\cA_I^*\ty_I^k}
	+ \inprod{X^k}{Z-\tZ^k}
	\\
	 && -\inprod{s^k}{y_I-\ty_I^k}
	  +  \frac{3}{2}\norm{\cA_E^*y_E-\cA_E^*\ty_E^k}^2+ \frac{3}{2}\norm{\cA_I^*y_I-\cA_I^*\ty_I^k}^2  + \frac{3}{2}\norm{Z-\tZ^k}^2  + \frac{1}{2}\norm{y_I-\ty_I^k}^2
	\\
	& :=  & \phi^k(y_E,y_I,Z),
\end{eqnarray*}
where $X^k = \Pi_{\S^n_+}(\cA_E^*\ty_E^k+\cA_I^*\ty_I^k+\tZ^k+G)$, $s^k = \Pi_{\cK}(g-\ty_I^k)$.

At each iteration of the APG method \cite{beck2009fast,tseng2008accelerated} applied to
($\widehat{\bD}$), we need to solve the following minimization
subproblem at the $k$th iteration:
\begin{eqnarray*}
	&\min_{y_E,y_I, Z} \Big\{ \phi^k (y_E,y_I,Z)\Big\}  &
\end{eqnarray*}
whose optimal solution is given by
\begin{eqnarray}\label{eq:eq-ineq}
& &\left\{
\begin{array}{l}		
y_E^{k} = (\cA_E \cA_E^*)^{-1} \big(\frac{1}{3}(b_E-\cA X^k)+\cA_E \cA_E^* \ty_E^k \big),  \\[5pt]	
y_I^{k} = (\cA_I \cA_I^*+\frac{1}{3}\cI_{m_I})^{-1} \big(\frac{1}{3}(s^k-\cA X^k+\ty_I^k)+\cA_I \cA_I^*\ty_I^k \big),  \\[5pt]		
Z^{k} = \Pi_\cP(\frac{1}{3}X^k-\tZ^k)-(\frac{1}{3}X^k-\tZ^k).
\end{array}\right.
\end{eqnarray}
With the computed iterate $(y_E^k,y_I^k,Z^k)$, we update the
auxiliary iterate $(\ty_E^{k+1},\ty_I^{k+1},\tZ^{k+1})$ similarly as in
Step 2 of Algorithm ABCD-1.

\subsection{An enhanced ARBCG method for solving (D)}

Next we describe the enhanced
accelerated randomized block coordinate gradient (ARBCG) method
  for solving (\bD). Our algorithm is modified from
Algorithm 3 in \cite{lin2014accelerated} for the sake of numerical comparison, and the steps
are given as follows.

\bigskip
\centerline{\fbox{\parbox{\textwidth}{
			{\bf Algorithm eARBCG}: {\bf An enhanced ARBCG method for (\bD).}
			\\[5pt]
We use the notation in \eqref{prob-gen} with $q=4$ and
$x= ((Z,v),S,y_E,y_I)$.
			Select an initial point  $ \tx^1
= x^1 $. Set $\alpha_0=1/q $ and
			$k=1$. Iterate the following steps.
			\begin{description}				
				\item[Step 1.]  Compute $\alpha_k = \frac{1}{2}\Big( \sqrt{\alpha_{k-1}^4 + 4\alpha_{k-1}^2}-\alpha_{k-1}^2\Big)$.
				\item [Step 2.] Compute
               	$\widehat{x}^{k+1} = (1-\alp_k)x^k +\alp_k \tx^k. $				
\item[Step 3.] Choose $i_k\in \{ 1,\dots,q\}$ uniformly at random and compute
\begin{eqnarray*}
 \tx^{k+1}_{i_k} = \mbox{argmin} \Big\{ \inprod{\grad_{i_k} \zeta(\widehat{x}^k)}{x_{i_k}-
\widehat{x}^k_{i_k}} + \frac{q\alp_k}{2} \norm{x_{i_k}-\tx^k_{i_k}}_{\cT_{i_k}}^2
+ \theta_{i_k}(x_{i_k})
\Big\},
\end{eqnarray*}
where $\cT_1$, $\cT_2$ are identity operators, and
$\cT_3 = \cA_E\cA_E^*$ and $\cT_4 = \cA_I\cA_I^*+\cI$.
Here $\theta_1(x_1) = \delta_\cP^*(-Z)+ \delta_\cK^*(-v)$,
$\theta_2(x_2) = \delta_{\S^n_+}(S)$,   $\theta_3(x_3) = 0 = \theta_4(x_4)$, and $\zeta$ is the smooth part of $F$ in (\bD).
Set   {$\tx^{k+1}_{i} = \tx^{k}_{i}$ for all $i \neq i_k$ and}
\begin{eqnarray*}
x^{k+1}_i = \left\{ \begin{array}{ll}
 \widehat{x}^k_{i} + q\alp_k (\tx^{k+1}_i - \tx^{k}_i) &\mbox{if $i=i_k$},
\\[5pt]
\widehat{x}^k_i & \mbox{if $i\not=i_k$}.
\end{array} \right.
\end{eqnarray*}
			\end{description}
}}}

\bigskip
Note that in the original ARBCG algorithm in \cite{lin2014accelerated},
the linear operators $\cT_3$ and $\cT_4$ are fixed to be
$\cT_3 = \lam_{\max}(\cA_E\cA_E^*)\cI$ and
$\cT_4 = \lam_{\max}(\cA_I\cA_I^* + \cI)\cI$. Our enhancement to the
algorithm is in using the operators
$\cT_3 = \cA_E\cA_E^*$ and $\cT_4 = \cA_I\cA_I^* + \cI$.
Indeed the practical performance of the eARBCG with the latter choices of $\cT_3$ and  $\cT_4$
is much better than the former more conservative choices. However,
we should note that although  the non-convergence  of the  eARBCG  method is never  observed
for the problems tested in our numerical experiments, the theoretical
convergence of the eARBCG has yet to be established,  for which we leave as a future research topic.

\section{Numerical experiments}
\label{sec-numerical-ABCD}

In our numerical experiments,
we test the algorithms designed  in the last two sections
to the  least squares semidefinite programming problem $\mbox{({\bf P})}$
by taking $G=-C, g=0$ for the data arising from various
classes of  SDP problems
of the form given in \eqref{eq-SDP}. Specifically, the LSSDP problem
corresponds to the first subproblem \eqref{eq-ppa} of the PPA for solving \eqref{eq-SDP} by
setting
$k=0$, $X^0=0$, $s^0=0$ and $\sig_0 = 1$.

\subsection{SDP problem sets}

Now we describe the classes of SDP problems
we considered in our numerical experiments.

\bigskip
(i) SDP problems coming from
the relaxation of a binary integer nonconvex quadratic (BIQ) programming:
\begin{eqnarray}
\min\Big\{ \frac{1}{2}x^TQx + \inprod{c}{x} \mid  x\in\{0,1\}^{n-1}\Big\}.
\end{eqnarray}
This problem has  been shown in \cite{burer2009copositive}  that under some mild
assumptions, it can equivalently be reformulated as the following
completely positive programming (CPP) problem:
\begin{eqnarray}
\min\Big\{ \frac{1}{2} \inprod{Q}{X_0} + \inprod{c}{x} \mid
\diag(X_0) = x, X = [X_0, x; x^T, 1]\in
{\cal C}_{pp}^n\Big\},
\end{eqnarray}
where
${\cal C}_{pp}^n$ denotes the $n$-dimensional completely positive cone.
It is well known that even though ${\cal C}_{pp}^n$ is convex, it is
computationally intractable. To solve the CPP problem,  one can
relax ${\cal C}_{pp}^n$ to $\S^n_+\cap \cN$, and the relaxed problem
has the form \eqref{eq-SDP}:
\begin{eqnarray}
\begin{array}{ll}
\min & \frac{1}{2} \inprod{Q}{X_0} + \inprod{c}{x} \\[5pt]
{\rm s.t.} &  {\rm diag}(X_0) - x = 0, \; \alpha = 1,
\quad X = \left[ \begin{array}{cc} X_0 & x \\ x^T & \alpha  \end{array}\right] \in \cS^n_+,\quad
X \in \cP,
\end{array}
\label{eq-BIQ-SDP}
\end{eqnarray}
where the polyhedral cone $\cP = \{X \in \cS^n \mid X \geq 0 \}$. In our numerical experiments, the test data for $Q$ and $c$ are taken from Biq Mac Library  maintained by Wiegele, which is available at {\tt http://biqmac.uni-klu.ac.at/biqmaclib.html}.

\bigskip
(ii) SDP problems arising from the relaxation of maximum stable set problems. Given a graph $G$ with edge set $\mathcal{E}$, the SDP relaxation  $\theta_+(G)$ of the maximum stable set problem is given by
\begin{eqnarray}
\theta_+(G) = \max \{ \inprod{ee^T}{X} \mid \;\inprod{E_{ij}}{X} =0,
(i,j) \in \mathcal{E}, \,\inprod{I}{X} =1, \, X \in \cS^n_+,\, X \in \cP \},\quad
\label{eq-theta-plus}
\end{eqnarray}
where $E_{ij}=e_ie_j^{\T} + e_je_i^{\T}$ and $e_i$ denotes the $i$th column
of the identity matrix, $\cP = \{X \in \cS^n \mid X \geq 0 \}$. In our numerical experiments, we test the graph instances $G$ considered in
\cite{Sloane}, \cite{toh2004}, and \cite{DIMACS}.

\bigskip
(iii) SDP relaxation for computing lower
bounds for quadratic assignment problems (QAPs).  Given
matrices $A, B \in \S^n$, the QAP is given by
$
v^*_{\rm QAP} := \min \{ \langle X, AXB \rangle \, :\, X \in \varPi \}
$
where $\varPi$ is the set of $n\times n$ permutation matrices.

For a matrix $X = [x_1,\dots,x_n]\in \Re^{n\times n}$, we will
identify it with the $n^2$-vector $x = [x_1; \dots; x_n]$. For a
matrix $Y\in R^{n^2 \times n^2}$, we let $Y^{ij}$ be the $n\times n$
block corresponding to $x_i x_j^T$ in the matrix $xx^T$.
In \cite{povh2009copositive}, it is shown that $v^*_{\rm QAP}$ is bounded below by the
following number generated from the SDP relaxation:
\begin{eqnarray}
\begin{array}{rl}
v \,:=\, \min & \inprod{B \otimes A}{Y} \\[5pt]
{\rm s.t.} &  \sum_{i=1}^n Y^{ii} = I, \,
\inprod{I}{Y^{ij}} = \delta_{ij} \quad \forall \, 1\leq i \leq j\leq n, \\[5pt]
&
\inprod{E}{Y^{ij}} = 1\quad \forall\, 1\leq i \leq j\leq n,
\quad Y \in \cS^n_+,\; Y \in \cP,
\end{array}
\label{eq-PR}
\end{eqnarray}
where the sign $``\otimes"$ stands for the Kronecker product, $E$ is the matrix of ones, and $\delta_{ij} = 1$ if $i=j$,
and $0$ otherwise, $\cP = \{X \in \cS^{n^2} \mid X \geq 0 \}$. In our numerical experiments, the test instances $(A,B)$ are taken
from the QAP Library \cite{QAPLIB}.

\bigskip
(iv) SDP relaxations of  clustering problems (RCPs)  described in \cite[eq. (13)]{peng2007approximating}:
\begin{eqnarray}
\min \Big\{ \inprod{-W}{X} \mid Xe = e, \inprod{I}{X} = K, X \in \S^n_+,
X \in \cP \Big\},
\end{eqnarray}
where $W$ is the so-called affinity matrix whose entries represent the pairwise similarities of the objects in the dataset, $e$ is the vector of ones, and $K$ is the number of clusters, $\cP = \{X \in \cS^n \mid X \geq 0 \}$.  All the data sets we tested are from the UCI Machine Learning Repository (available at {\tt http://archive.ics.uci.edu/ml/datasets.html}). For some large data instances, we only select the first $n$ rows. For example, the original data instance ``{\tt spambase}" has 4601 rows, we select the first 1500 rows to obtain the test problem ``{\tt spambase-large.2}"
for which the number ``{\tt 2}" means that there are $K = 2$ clusters.

\bigskip
(v) SDP problems arising from the SDP relaxation of frequency assignment problems (FAPs) \cite{eisenblatter2002frequency}.  Given a network represented by a graph $G$ with edge set $\cE$ and an edge-weight matrix $W$, a certain type of frequency assignment problem on $G$  can be relaxed into the following SDP (see \cite[eq. (5)]{burer2003computational}):
\begin{eqnarray} \label{eq-fap-sdpp}
\begin{array}{ll}
\max & \inprod{(\frac{k-1}{2k})L(W)-\frac{1}{2}{\rm Diag}(We)}{ X}
\\[5pt]
{\rm s.t.} & {\rm diag}(X) = e, \;\; X  \in \S^n_+,\;\; X \in \cP,
\end{array}
\end{eqnarray}
where  $L(W) := {\rm Diag}(We)-W$ is the Laplacian matrix,  $e \in \Re^n$ is the vector of all ones, and
 $\cP = \{X \in \cS^n \mid L \leq X \leq U \}$,
\begin{eqnarray*}
	L_{ij} =
	\begin{cases}
		-\frac{1}{k-1} &  \forall (i,j) \in \cE, \\
		- \infty & {\rm otherwise},
	\end{cases}
	\quad
	U_{ij} =
	\begin{cases}
		-\frac{1}{k-1} &  \forall (i,j) \in \cU, \\
		\infty & {\rm otherwise}
	\end{cases}
\end{eqnarray*}
with $k>1$ being a given integer and $\cU$ is a given subset of $\cE$.

\bigskip
(vi) For the SDP problems described in (\ref{eq-BIQ-SDP}) arising from relaxing the BIQ problems, in order to get tighter bounds, we may add  in some valid inequalities to  get the following problems:
\begin{eqnarray}\begin{array}{ll}
\min & \frac{1}{2} \inprod{Q}{Y} + \inprod{c}{x} \\[5pt]
{\rm s.t.} &  {\rm diag}(Y) - x = 0, \; \alpha = 1, \;
X = \left[ \begin{array}{cc} Y & x \\ x^T & \alpha  \end{array}\right]\in \S^n_+,\;
X \in \cP,
\\[8pt]
& 0 \leq -Y_{ij} + x_i \leq 1,\;
0 \leq -Y_{ij} + x_j \leq 1,\;
-1 \leq Y_{ij} - x_i -x_j \leq 0,\;
\\[8pt]
& \forall \; 1\leq i < j, \; j\leq n-1,
\end{array}
\label{eq-BIQineq-SDP}
\end{eqnarray}
where $\cP = \{X \in \cS^n \mid X \geq 0 \}$. For convenience, we call the problem in (\ref{eq-BIQineq-SDP}) an extended BIQ problem.
Note that the last set of inequalities in \eqref{eq-BIQineq-SDP}
are obtained from the valid inequalities $0 \leq x_i(1-x_j) \leq 1$, $0 \leq x_j(1-x_i) \leq 1$, $0 \leq (1-x_i)(1-x_j) \leq 1$
when $x_i,x_j$ are binary variables.

\subsection{Numerical results}

In this section, we compare the performance of the ABCD, APG, eARBCG and BCD methods
for solving the LSSDP (\bP). Note that the BCD method follows the template
described in \eqref{bcd}  {with $q=4$ and $x= ((Z,v),S,y_E,y_I)$}.
All our computational results  are obtained by running {\sc Matlab} on a  workstation (20-core, Intel Xeon E5-2670 v2 @ 2.5GHz, 64GB RAM).

In our numerical experiments, for problem (\bP), we assume that $\gamma :=\max \{\norm{G}, \norm{g} \} \leq 1$, otherwise we can solve an equivalent rescaled problem
\begin{eqnarray*}
	\mbox{($\overline{\bf {P}}$)} & &
	\begin{array}{ll}	
		\min &  \frac{1}{2}\norm{X-\bar{G}}^2 + \frac{1}{2}\norm{s-\bar{g}}^2  \\[5pt]
		{\rm s.t.} &  \cA_E(X) = \bar{b}_E, \; \cA_I X -s = 0, \; X\in \S^n_+, \; X \in \bar{\cP}, \; s \in \bar{\cK},
	\end{array}
\end{eqnarray*}
where $\bar{G} = \frac{G}{\gamma}$, $\bar{g} = \frac{g}{\gamma}$,   $\bar{b}_E = \frac{b_E}{\gamma}$, $\bar{\cP} = \{X \mid \gamma X \in \cP \}$, $\bar{\cK} = \{s \mid \gamma s \in \cK \}$. Note that $(X,s)$ is a solution to (\bP) if and only if $(\frac{X}{\gamma}, \frac{s}{\gamma})$ is a solution to ($\overline{\bf{P}}$).

Note that under a suitable Slater's condition,
the KKT conditions for (\bP) and (\bD) are given as follows:
\begin{eqnarray}
\begin{array}{rcl}
A_E X &=& b_E,
 \;\;  \cA_I X - s=  0, \;\; X-Y = 0, \\[5pt]
X & = & \Pi_{\S^n_+}(\cA_E^*y_E+\cA_I^* y _I+Z+G),  \\[5pt]
Y &=& \Pi_{\cP}(\cA_E^*y_E+\cA_I^* y _I+S+G),  \\[5pt]
s &=& \Pi_{\cK}(g-y_I).
\end{array}
\label{eq-KKT}
\end{eqnarray}
Thus
we measure the accuracy of an approximate optimal solution $(Z,v,S,y_E,y_I)$ for
({\bf{D}}) by using the following relative residual:
\begin{eqnarray}
	\eta = \max \{\eta_1, \eta_2, \eta_3\},  \label{eq-eta}
\end{eqnarray}
where $\eta_1 =\frac{\norm{b_E-\cA_E X}}{1+\norm{b_E}}$,
$\eta_2 =\frac{\norm{X-Y}}{1+\norm{X}}$, $\eta_3 =\frac{\norm{s-\cA_I X}}{1+\norm{s}}$, $X = \Pi_{\S^n_+}(\cA_E^*y_E+\cA_I^* y _I+Z+G)$,
$Y = \Pi_{\cP}(\cA_E^*y_E+\cA_I^* y _I+S+G)$,
$s = \Pi_{\cK}(g-y_I)$.
Additionally, we compute the relative gap defined by
\begin{eqnarray}
	&\eta_g = \frac{\frac{1}{2}\norm{X-G}^2 + \frac{1}{2}\norm{s-g}^2 - F(Z,v,S,y_E,y_I) }{1+\frac{1}{2}\norm{X-G}^2 + \frac{1}{2}\norm{s-g}^2 + \abs{F(Z,v,S,y_E,y_I)}}.&
	\label{eq-gap}
\end{eqnarray}
Let  $\varepsilon > 0$ be a given  accuracy tolerance. We terminate the
ABCD, APG, eARBCG and BCD methods when
$\eta  < \varepsilon. $

Table \ref{table-numsolved} shows the number  of problems that  have been successfully solved to the accuracy of  $10^{-6}$ in $\eta$ by each of  the four solvers: ABCD, APG, eARBCG and BCD, with the maximum number of iterations set at $25000$\footnote{The maximum number of iterations for eARBCG is $25000 q$, where $q$ is the number of block variables, i.e., $q = 4$ for extended BIQ problems and $q = 3$ otherwise, since eARBCG only updates one block variable at each iteration.}. As can be seen, only ABCD  can solve all the problems to the desired accuracy of  $10^{-6}$. The performance of the BCD method is especially poor, and it can only solve $201$ problems out of $616$ to the desired accuracy.

Table \ref{table-numsolved_tol} shows the number  of problems  that  have been successfully solved to the accuracy of  $10^{-6}$, $10^{-7}$ and $10^{-8}$  in $\eta$ by the solver ABCD, with the maximum number of iterations set at $25000$. As can be seen,  ABCD can even solve almost all the problems to the high accuracy of $10^{-8}$.

\begin{table}[h]
	\begin{center}
		\caption{Number of problems which are solved to the accuracy of  $ 10^{-6}$ in $\eta$.} \label{table-numsolved}
		\begin{tabular}{| c | c | c | c | c |}
			\hline
			\multicolumn{1}{|c}{problem set (No.) $\backslash$ solver} & \multicolumn{1}{|c|}{ABCD} & \multicolumn{1}{|c|}{APG} & \multicolumn{1}{|c|}{eARBCG}  & \multicolumn{1}{|c|}{BCD} \\ \hline	
			$\theta_{+}$ (64) & 64 & 64 & 64 & 11 \\ \hline FAP ( 7) &  7 &  7 &  7 &  7 \\ \hline QAP (95) & 95 & 95 & 24 &  0 \\ \hline BIQ (165) & 165 & 165 & 165 & 65 \\ \hline RCP (120) & 120 & 120 & 120 & 108 \\ \hline exBIQ (165) & 165 & 141 & 165 & 10 \\ \hline Total (616) & 616 & 592 & 545 & 201 \\ \hline 
		\end{tabular}
	\end{center}
\end{table}

\begin{table}[h]
	\begin{center}
		\caption{Number of problems which are solved to the accuracy of  $ 10^{-6}$, $10^{-7}$ and $10^{-8}$, in $\eta$ by the ABCD method.} \label{table-numsolved_tol}
		\begin{tabular}{| c | c | c | c | c |}
			\hline
			\multicolumn{1}{|c}{problem set (No.) $\backslash$ $\varepsilon$} & \multicolumn{1}{|c|}{$10^{-6}$} & \multicolumn{1}{|c|}{$10^{-7}$} & \multicolumn{1}{|c|}{$10^{-8}$}   \\ \hline
			$\theta_{+}$ (64) & 64 & 58 & 52  \\ \hline FAP ( 7) &  7 &  7 &  7  \\ \hline QAP (95) & 95 & 95 & 95 \\ \hline BIQ (165) & 165 & 165 & 165  \\ \hline RCP (120) & 120 & 120 & 118  \\ \hline exBIQ (165) & 165 & 165 & 165  \\ \hline Total (616) & 616 & 610 & 602  \\ \hline 
		\end{tabular}
	\end{center}
\end{table}
Table \ref{table_ABCD_APG_ARBCG}  compares the performance of the
solvers  ABCD, APG and eARBCG on a subset of the 616 tested
LSSDP  problems using the tolerance $\varepsilon = 10^{-6}$.
The full table for all the 616 problems is available at
\url{http://www.math.nus.edu.sg/~mattohkc/publist.html}. 
The first three columns of Table \ref{table_ABCD_APG_ARBCG} give the problem name, the dimensions of the variables $y_E$ ($m_E$) and $y_I$ ($m_I$), the size of the matrix $C$ ($n_s$)  in (P), respectively.  The number of iterations\footnote{Note that the two numbers of iterations in ABCD are for  ABCD-2 and ABCD-1, respectively.}, the relative residual $\eta$ and relative gap $\eta_g$, as well as computation times (in the format hours:minutes:seconds) are listed in the last twelve columns. As can be seen, ABCD is much faster than APG and eARBCG for most of the problems.

\begin{figure}[h]
	\begin{center}
		\ifpdf
		\includegraphics[width=0.5\textwidth]{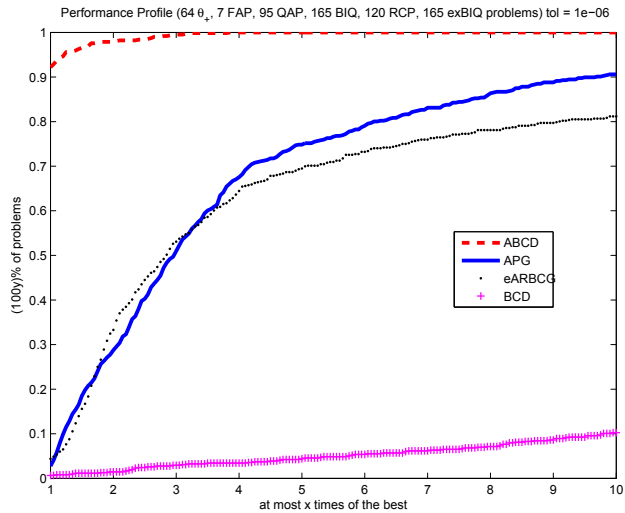}
		\else
		\includegraphics[width=0.5\textwidth]{Performance_Profile.eps}
		\fi
		\caption{Performance profiles of ABCD, APG, eARBCG and BCD on $[1,10].$}
		\label{fig-all}
	\end{center}
\end{figure}

Figure \ref{fig-all} shows the performance profiles of the ABCD, APG, eARBCG and BCD methods for {all  the $616$ tested problems}. We recall that a point $(x,y)$ is in the performance profile curve of a method if and only if it can solve exactly $(100y)\%$ of all the tested problems at most $x$ times slower than any other methods. It can be seen that ABCD outperforms the other $3$ methods by a significant margin. Furthermore, the ABCD method is more than ten times faster than the BCD method for vast majority of the  problems. It
is quite surprising that a simple novel acceleration step with a special BCD cycle can improve the
performance of the standard Gauss-Seidel BCD method by such a
dramatic margin.

\begin{figure}[h]
	\begin{center}
		\ifpdf
		\includegraphics[width=0.5\textwidth]{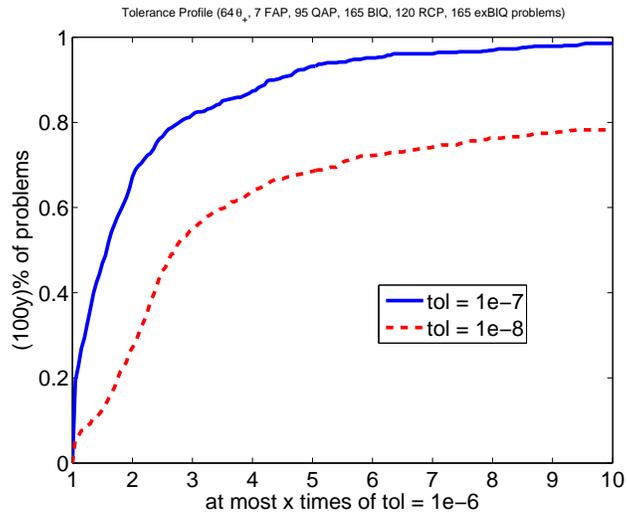}
		\else
		\includegraphics[width=0.5\textwidth]{Tolerance_Profile.eps}
		\fi
		\caption{Tolerance profiles of ABCD on $[1,10].$}
		\label{fig-tolerance}
	\end{center}
\end{figure}

Figure \ref{fig-tolerance} shows the tolerance  profiles of the ABCD method
for {all  the $616$ tested problems}. Note that a point $(x,y)$ is in the tolerance profile curve if and only if it can solve exactly $(100y)\%$ of all the tested problems at most $x$ times slower than the time taken to reach the tolerance of $10^{-6}$.

\section{Conclusions}
\label{sec-conclusion-ABCD}

We have designed an inexact accelerated block coordinate gradient descent (ABCGD)
method for solving a multi-block convex minimization problem
whose objective is the sum of a coupled smooth function with Lipschitz continuous
gradient and a separable (possibly nonsmooth) function involving
only the first two blocks.
An important class of problems with the specified structure
is the dual of
least squares semidefinite programming (LSSDP)
where the primal matrix variable must satisfy given linear equality and inequality
constraints, and  must also lie in the intersection of the
cone of positive semidefinite matrices and a simple polyhedral set.

Our inexact  ABCGD method has $O(1/k^2)$ iteration complexity
if the subproblems are solved progressively more accurately.
The design of our    ABCGD method relied on recent advances
in the symmetric Gauss-Seidel technique for solving a
convex minimization problem whose objective is the sum of a
 multi-block quadratic function  and a
non-smooth function involving only the first block.
Extensive numerical experiments on various classes of over 600 large scale
LSSDP problems demonstrate that our ABCGD method, which reduces to the ABCD method for LSSDP problems,    not only
can solve the problems to high accuracy, but it is also far more efficient
than (a) the well known BCD method,
(b) the eARBCG (an enhanced version of the accelerated randomized block coordinate gradient) method,
and (c) the APG (accelerated proximal gradient) method.

\section*{Acknowledgements}
We  would like to thank Ms Ying Cui from National University of Singapore for her  comments on refining the Danskin-type Theorem presented  in a previous version of this paper.


\bibliographystyle{siam}
\bibliography{LSSDPbib}

\begin{thebibliography}{10}

\bibitem{beck2009fast}
{\sc A.~Beck and M.~Teboulle}, {\em A fast iterative shrinkage-thresholding
  algorithm for linear inverse problems}, SIAM Journal on Imaging Sciences, 2
  (2009), pp.~183--202.

\bibitem{beck2013convergence}
{\sc A.~Beck and L.~Tetruashvili}, {\em On the convergence of block coordinate
  descent type methods}, SIAM Journal on Optimization, 23 (2013),
  pp.~2037--2060.

\bibitem{burer2009copositive}
{\sc S.~Burer}, {\em On the copositive representation of binary and continuous
  nonconvex quadratic programs}, Mathematical Programming, 120 (2009),
  pp.~479--495.

\bibitem{burer2003computational}
{\sc S.~Burer, R.~D. Monteiro, and Y.~Zhang}, {\em A computational study of a
  gradient-based log-barrier algorithm for a class of large-scale \mbox{SDPs}},
  Mathematical Programming, 95 (2003), pp.~359--379.

\bibitem{chambolle2015remark}
{\sc A.~Chambolle and T.~Pock}, {\em A remark on accelerated block coordinate
  descent for computing the proximity operators of a sum of convex functions},
  Optimization Online,  (2015).

\bibitem{Danskin66}
{\sc J.~M. Danskin}, {\em The theory of min-max, with applications}, SIAM
  Journal on Applied Mathematics, 14 (1966), pp.~641--664.

\bibitem{eisenblatter2002frequency}
{\sc A.~Eisenbl{\"a}tter, M.~Gr{\"o}tschel, and A.~M. Koster}, {\em Frequency
  planning and ramifications of coloring}, Discussiones Mathematicae Graph
  Theory, 22 (2002), pp.~51--88.

\bibitem{facchinei2003finite}
{\sc F.~Facchinei and J.-S. Pang}, {\em Finite-dimensional variational
  inequalities and complementarity problems}, vol.~1, Springer
  Science$+$Business Media, 2003.

\bibitem{fercoq2013accelerated}
{\sc O.~Fercoq and P.~Richt{\'a}rik}, {\em Accelerated, parallel and proximal
  coordinate descent}, arXiv preprint arXiv:1312.5799,  (2013).

\bibitem{grippo2000convergence}
{\sc L.~Grippo and M.~Sciandrone}, {\em On the convergence of the block
  nonlinear \mbox{Gauss--Seidel} method under convex constraints}, Operations
  Research Letters, 26 (2000), pp.~127--136.

\bibitem{QAPLIB}
{\sc P.~Hahn and M.~Anjos}, {\em {QAPLIB} -- a quadratic assignment problem
  library}.
\newblock \url{http://www.seas.upenn.edu/qaplib}.

\bibitem{hiriart1984generalized}
{\sc J.-B. Hiriart-Urruty, J.-J. Strodiot, and V.~H. Nguyen}, {\em Generalized
  {Hessian} matrix and second-order optimality conditions for problems with
  {$C^{1,1}$} data}, Applied Mathematics and Optimization, 11 (1984),
  pp.~43--56.

\bibitem{jiang2012inexact}
{\sc K.~Jiang, D.~F. Sun, and K.-C. Toh}, {\em An inexact accelerated proximal
  gradient method for large scale linearly constrained convex \mbox{SDP}}, SIAM
  Journal on Optimization, 22 (2012), pp.~1042--1064.

\bibitem{aug2QP}
{\sc X.~Li, D.~F. Sun, and K.-C. Toh}, {\em \mbox{QP-PAL}: A 2-phase proximal
  augmented {Lagrangian} method for high dimensional convex quadratic
  programming problems}, preprint,  (2015).

\bibitem{LiSunToh2015}
\leavevmode\vrule height 2pt depth -1.6pt width 23pt, {\em A \mbox{S}chur
  complement based semi-proximal \mbox{ADMM} for convex quadratic conic
  programming and extensions}, Mathematical Programming, to appear (2015).

\bibitem{lin2014accelerated}
{\sc Q.~Lin, Z.~Lu, and L.~Xiao}, {\em An accelerated proximal coordinate
  gradient method and its application to regularized empirical risk
  minimization}, arXiv preprint arXiv:1407.1296,  (2014).

\bibitem{nesterov1983method}
{\sc Y.~Nesterov}, {\em A method of solving a convex programming problem with
  convergence rate {$O(1/k^2)$}}, in Soviet Mathematics Doklady, vol.~27, 1983,
  pp.~372--376.

\bibitem{nesterov2004introductory}
\leavevmode\vrule height 2pt depth -1.6pt width 23pt, {\em Introductory
  lectures on convex optimization: a basic course}, vol.~87, Springer
  Science$+$Business Media, 2004.

\bibitem{nesterov2005smooth}
\leavevmode\vrule height 2pt depth -1.6pt width 23pt, {\em Smooth minimization
  of non-smooth functions}, Mathematical Programming, 103 (2005), pp.~127--152.

\bibitem{nesterov2012efficiency}
{\sc Y.~Nesterov}, {\em Efficiency of coordinate descent methods on huge-scale
  optimization problems}, SIAM Journal on Optimization, 22 (2012),
  pp.~341--362.

\bibitem{peng2007approximating}
{\sc J.~Peng and Y.~Wei}, {\em Approximating k-means-type clustering via
  semidefinite programming}, SIAM Journal on Optimization, 18 (2007),
  pp.~186--205.

\bibitem{povh2009copositive}
{\sc J.~Povh and F.~Rendl}, {\em Copositive and semidefinite relaxations of the
  quadratic assignment problem}, Discrete Optimization, 6 (2009), pp.~231--241.

\bibitem{rockafellar1970convex}
{\sc R.~T. Rockafellar}, {\em Convex analysis}, Princeton Mathematical Series,
  No. 28, Princeton University Press, Princeton, N.J., 1970.

\bibitem{rockafellar1974conjugate}
{\sc R.~T. Rockafellar}, {\em Conjugate duality and optimization}, vol.~14,
  SIAM, 1974.

\bibitem{rockafellar1976augmented}
\leavevmode\vrule height 2pt depth -1.6pt width 23pt, {\em Augmented
  \mbox{Lagrangians} and applications of the proximal point algorithm in convex
  programming}, Mathematics of Operations Research, 1 (1976), pp.~97--116.

\bibitem{rockafellar1976monotone}
\leavevmode\vrule height 2pt depth -1.6pt width 23pt, {\em Monotone operators
  and the proximal point algorithm}, SIAM Journal on Control and Optimization,
  14 (1976), pp.~877--898.

\bibitem{saha2013nonasymptotic}
{\sc A.~Saha and A.~Tewari}, {\em On the nonasymptotic convergence of cyclic
  coordinate descent methods}, SIAM Journal on Optimization, 23 (2013),
  pp.~576--601.

\bibitem{sardy2000block}
{\sc S.~Sardy, A.~G. Bruce, and P.~Tseng}, {\em Block coordinate relaxation
  methods for nonparametric wavelet denoising}, Journal of Computational and
  Graphical Statistics, 9 (2000), pp.~361--379.

\bibitem{Sloane}
{\sc N.~Sloane}, {\em Challenge problems: Independent sets in graphs}.
\newblock \url{http://www.research.att.com/∼njas/doc/graphs.html}, 2005.

\bibitem{toh2004}
{\sc K.-C. Toh}, {\em Solving large scale semidefinite programs via an
  iterative solver on the augmented systems}, SIAM Journal on Optimization, 14
  (2004), pp.~670--698.

\bibitem{DIMACS}
{\sc M.~Trick, V.~Chvatal, B.~Cook, D.~Johnson, C.~McGeoch, and R.~Tarjan},
  {\em The second dimacs implementation challenge --- \mbox{NP} hard problems:
  Maximum clique, graph coloring, and satisfiability}.
\newblock \url{http://dimacs.rutgers.edu/Challenges/}, 1992.

\bibitem{tseng1993dual}
{\sc P.~Tseng}, {\em Dual coordinate ascent methods for non-strictly convex
  minimization}, Mathematical Programming, 59 (1993), pp.~231--247.

\bibitem{tseng2001convergence}
\leavevmode\vrule height 2pt depth -1.6pt width 23pt, {\em Convergence of a
  block coordinate descent method for nondifferentiable minimization}, Journal
  of Optimization Theory and Applications, 109 (2001), pp.~475--494.

\bibitem{tseng2008accelerated}
\leavevmode\vrule height 2pt depth -1.6pt width 23pt, {\em On accelerated
  proximal gradient methods for convex-concave optimization}, manuscript,
  (2008).

\bibitem{tseng2009coordinate}
{\sc P.~Tseng and S.~Yun}, {\em A coordinate gradient descent method for
  nonsmooth separable minimization}, Mathematical Programming, 117 (2009),
  pp.~387--423.

\bibitem{yang2014sdpnal}
{\sc L.-Q. Yang, D.~F. Sun, and K.-C. Toh}, {\em {SDPNAL$+$}: A majorized
  semismooth \mbox{Newton-CG} augmented \mbox{Lagrangian} method for
  semidefinite programming with nonnegative constraints}, Mathematical
  Programming Computation,  (to appear).

\bibitem{zhao2010newton}
{\sc X.~Y. Zhao, D.~F. Sun, and K.-C. Toh}, {\em A \mbox{Newton-CG} augmented
  \mbox{Lagrangian} method for semidefinite programming}, SIAM Journal on
  Optimization, 20 (2010), pp.~1737--1765.

\end{thebibliography}

\newpage

\begin{landscape}
\begin{center}	
	\begin{tiny}
		\begin{longtable}{| ccc | ccc | ccc | ccc | ccc |}
			\caption{Performance of ABCD, APG and eARBCG on $\theta_+$, FAP, QAP, BIQ, RCP and extended BIQ problems ($\varepsilon = 10^{-6}$)} \label{table_ABCD_APG_ARBCG} \\
			\hline
			\multicolumn{3}{|c}{} & \multicolumn{3}{|c|}{iteration} & \multicolumn{3}{|c|}{$\eta$} & \multicolumn{3}{|c|}{$\eta_g$}  & \multicolumn{3}{|c|}{time} \\ \hline
			problem & $m_E;m_I$ & $n_{s}$ & ABCD & APG & eARBCG & ABCD & APG & eARBCG & ABCD & APG & eARBCG & ABCD & APG & eARBCG \\ \hline		
			\endhead
			 theta8 &  7905;0 &   400 &   0;1689 &  2583 &  7225 &  9.8-7  &  9.9-7  &  9.9-7  &  -1.3-7  &  -1.9-9  &  -3.2-7  & 3:24 & 5:20 & 11:02  \\ \hline 
 theta82 &  23872;0 &   400 &   0;1405 &  2227 &  6302 &  9.9-7  &  9.9-7  &  9.7-7  &  -1.5-7  &  4.3-9  &  -1.4-7  & 2:49 & 4:52 & 9:16  \\ \hline 
 theta83 &  39862;0 &   400 &   0;1672 &  2157 &  6242 &  9.9-7  &  9.9-7  &  9.9-7  &  1.6-8  &  5.8-9  &  -1.2-7  & 4:00 & 4:36 & 9:30  \\ \hline 
 theta10 &  12470;0 &   500 &   0;1820 &  2629 &  7839 &  9.9-7  &  9.9-7  &  9.9-7  &  -9.2-8  &  1.9-10  &  -2.5-7  & 7:19 & 9:08 & 18:26  \\ \hline 
 theta102 &  37467;0 &   500 &   0;1705 &  2266 &  6923 &  9.9-7  &  9.9-7  &  9.9-7  &  6.8-8  &  1.9-8  &  -1.5-7  & 7:04 & 8:19 & 16:25  \\ \hline 
 theta103 &  62516;0 &   500 &   0;1641 &  2228 &  6729 &  9.9-7  &  9.9-7  &  9.8-7  &  5.9-8  &  1.8-8  &  -1.3-7  & 6:50 & 8:31 & 17:01  \\ \hline 
 theta104 &  87245;0 &   500 &  25;1483 &  2235 &  6726 &  9.9-7  &  9.9-7  &  9.9-7  &  -1.8-7  &  1.6-8  &  -1.1-7  & 7:16 & 7:41 & 16:43  \\ \hline 
 theta12 &  17979;0 &   600 &   0;1961 &  2644 &  8728 &  9.9-7  &  9.8-7  &  9.7-7  &  -1.4-7  &  5.8-8  &  -2.6-7  & 12:28 & 14:40 & 27:32  \\ \hline 
 theta123 &  90020;0 &   600 &   0;1743 &  2366 &  7160 &  9.9-7  &  9.9-7  &  9.9-7  &  4.4-8  &  2.0-8  &  -1.3-7  & 11:27 & 13:37 & 27:35  \\ \hline 
 theta162 &  127600;0 &   800 &   0;2149 &  2670 &  7508 &  9.9-7  &  9.8-7  &  9.9-7  &  4.0-8  &  2.3-8  &  -1.3-7  & 29:59 & 30:58 & 54:54  \\ \hline 
 hamming-9-8 &  2305;0 &   512 &   0;382 &  615 &  6304 &  2.5-7  &  5.3-7  &  8.4-7  &  {\green{ -1.3-6}}  &  -2.6-7  &  -3.0-7  & 57 & 1:23 & 13:02  \\ \hline 
 hamming-10-2 &  23041;0 &  1024 &  89;254 &  817 &  5211 &  8.2-7  &  7.6-7  &  6.4-7  &  {\green{ 4.2-6}}  &  7.5-8  &  -5.4-8  & 4:39 & 6:28 & 56:51  \\ \hline 
 hamming-9-5-6 &  53761;0 &   512 &   0;253 &  456 &  1629 &  8.0-7  &  1.9-7  &  3.6-7  &  5.6-7  &  -1.2-8  &  -5.1-8  & 27 & 37 & 2:40  \\ \hline 
 brock400-1 &  20078;0 &   400 &   0;1586 &  2347 &  6553 &  9.9-7  &  9.9-7  &  9.6-7  &  -6.0-8  &  -3.4-9  &  -1.9-7  & 1:35 & 2:17 & 6:21  \\ \hline 
 G43 &  9991;0 &  1000 &   0;3350 &  3949 &  15965 &  9.9-7  &  9.9-7  &  9.9-7  &  -5.2-7  &  -4.6-8  &  -2.8-7  & 26:40 & 30:07 &   2:11:03  \\ \hline 
 G44 &  9991;0 &  1000 &   0;3351 &  4039 &  15965 &  9.9-7  &  9.8-7  &  9.9-7  &  -4.1-7  &  6.1-8  &  -2.9-7  & 28:02 & 31:43 &   2:34:37  \\ \hline 
 G45 &  9991;0 &  1000 &   0;3220 &  3868 &  15678 &  9.9-7  &  9.9-7  &  9.9-7  &  -4.6-7  &  -1.1-7  &  -2.0-7  & 25:28 & 29:15 &   2:51:51  \\ \hline 
 G46 &  9991;0 &  1000 &   0;3144 &  3647 &  15630 &  9.9-7  &  9.7-7  &  9.8-7  &  -4.7-7  &  1.9-7  &  -1.9-7  & 24:59 & 27:29 &   2:55:29  \\ \hline 
 G47 &  9991;0 &  1000 &   0;3378 &  4025 &  15607 &  9.9-7  &  9.8-7  &  9.9-7  &  -5.5-7  &  7.8-9  &  -2.3-7  & 26:55 & 32:19 &   1:28:28  \\ \hline 
 G51 &  5910;0 &  1000 &   0;11273 &  16058 &  41754 &  8.6-7  &  9.9-7  &  9.9-7  &  -5.2-7  &  -1.9-8  &  -9.1-7  &   1:39:35 &   2:18:14 &   4:03:33  \\ \hline 
 G52 &  5917;0 &  1000 &   0;7955 &  18194 &  36594 &  9.9-7  &  9.8-7  &  9.9-7  &  -5.0-7  &  7.3-9  &  -9.3-7  &   1:09:04 &   3:03:03 &   3:43:49  \\ \hline 
 G53 &  5915;0 &  1000 &   0;10887 &  16420 &  39485 &  9.6-7  &  9.5-7  &  9.8-7  &  -5.6-7  &  -1.6-8  &  -7.3-7  &   1:41:55 &   3:48:26 &   4:55:22  \\ \hline 
 G54 &  5917;0 &  1000 &   0;10883 &  13798 &  41754 &  9.5-7  &  9.9-7  &  9.9-7  &  -4.9-7  &  1.5-8  &  -9.2-7  &   2:46:39 &   4:09:19 &   7:36:36  \\ \hline 
 1dc.512 &  9728;0 &   512 &   0;3755 &  5543 &  17564 &  9.9-7  &  9.5-7  &  9.8-7  &  -3.6-7  &  -1.6-8  &  -3.9-7  & 6:21 & 17:53 & 25:12  \\ \hline 
 1et.512 &  4033;0 &   512 &   0;3265 &  3181 &  11786 &  9.9-7  &  9.9-7  &  9.9-7  &  -3.8-7  &  -1.8-8  &  -5.0-7  & 5:33 & 11:14 & 16:14  \\ \hline 
 1tc.512 &  3265;0 &   512 &   0;5590 &  10046 &  26847 &  9.9-7  &  9.8-7  &  9.8-7  &  -6.8-7  &  1.3-8  &  -4.0-7  & 9:39 & 33:42 & 37:02  \\ \hline 
 2dc.512 &  54896;0 &   512 &  16;2871 &  3468 &  8987 &  9.5-7  &  9.7-7  &  9.8-7  &  -3.8-7  &  7.2-9  &  -2.4-7  & 5:24 & 12:16 & 12:47  \\ \hline 
 1zc.512 &  6913;0 &   512 &   0;3370 &  2666 &  18104 &  9.9-7  &  9.1-7  &  9.4-7  &  -3.2-7  &  -8.9-7  &  -3.8-7  & 4:59 & 7:27 & 24:15  \\ \hline 
 1dc.1024 &  24064;0 &  1024 &   0;3074 &  5976 &  16655 &  9.9-7  &  9.5-7  &  9.9-7  &  -2.7-7  &  3.4-8  &  -2.6-7  & 28:57 &   2:17:15 &   1:53:08  \\ \hline 
 1et.1024 &  9601;0 &  1024 &   0;5119 &  7874 &  26394 &  9.9-7  &  9.9-7  &  9.9-7  &  -3.7-7  &  -7.3-8  &  4.4-7  & 46:55 &   2:55:02 &   3:42:44  \\ \hline 
 1tc.1024 &  7937;0 &  1024 &   0;6282 &  10496 &  26839 &  9.9-7  &  9.9-7  &  9.9-7  &  -6.4-7  &  -9.6-10  &  -1.8-7  &   1:00:21 &   3:46:21 &   5:03:53  \\ \hline 
 1zc.1024 &  16641;0 &  1024 &   0;2281 &  1724 &  13029 &  9.7-7  &  9.9-7  &  9.7-7  &  {\green{ -1.3-6}}  &  -3.5-7  &  6.2-8  & 18:16 & 32:11 &   2:29:37  \\ \hline 
 2dc.1024 &  169163;0 &  1024 &  18;2029 &  4910 &  11689 &  9.9-7  &  9.7-7  &  9.7-7  &  -2.9-7  &  2.1-9  &  -3.2-7  & 19:23 &   1:46:34 & 46:54  \\ \hline 
 1dc.2048 &  58368;0 &  2048 &   0;2700 &  5447 &  17472 &  9.9-7  &  9.9-7  &  9.9-7  &  -8.5-8  &  -1.6-7  &  6.7-7  &   2:54:31 &   11:05:59 &   7:38:17  \\ \hline 
 1et.2048 &  22529;0 &  2048 &   0;5267 &  14243 &  30322 &  9.9-7  &  9.3-7  &  9.7-7  &  -4.2-7  &  -1.9-8  &  9.7-7  &   5:13:18 &   21:40:08 &   14:08:56  \\ \hline 
 1tc.2048 &  18945;0 &  2048 &   0;7737 &  15235 &  65261 &  9.8-7  &  9.8-7  &  9.4-7  &  -5.4-7  &  -5.1-8  &  {\green{ 1.4-6}}  &   7:34:43 &   22:17:52 &   31:38:14  \\ \hline 
 2dc.2048 &  504452;0 &  2048 &   0;5392 &  6996 &  14000 &  9.9-7  &  9.8-7  &  9.9-7  &  8.8-9  &  1.6-9  &  -4.3-7  &   5:38:09 &   6:15:33 &   6:56:29  \\ \hline 
 1zc.2048 &  39425;0 &  2048 &   0;2928 &  1730 &  32486 &  9.2-7  &  9.4-7  &  8.7-7  &  -9.3-8  &  -6.0-7  &  6.6-7  &   4:16:36 &   1:38:33 &   21:51:00  \\ \hline 
 fap08 &  120;0 &   120 &   0;47 &  111 &  223 &  8.4-7  &  9.5-7  &  9.7-7  &  -6.5-8  &  2.0-7  &  {\green{ 2.6-6}}  & 00 & 01 & 01  \\ \hline 
 fap09 &  174;0 &   174 &   0;54 &  111 &  258 &  9.4-7  &  9.6-7  &  9.1-7  &  1.4-7  &  1.9-7  &  {\green{ 2.1-6}}  & 01 & 01 & 02  \\ \hline 
 fap10 &  183;0 &   183 &   0;45 &  67 &  254 &  8.4-7  &  7.9-7  &  6.2-7  &  7.5-8  &  1.2-7  &  3.2-7  & 01 & 01 & 02  \\ \hline 
 fap11 &  252;0 &   252 &   0;53 &  76 &  219 &  8.6-7  &  9.9-7  &  9.2-7  &  1.8-7  &  2.1-7  &  {\green{ 1.1-6}}  & 01 & 01 & 04  \\ \hline 
 fap12 &  369;0 &   369 &   0;58 &  94 &  258 &  7.8-7  &  9.7-7  &  7.1-7  &  -9.2-8  &  4.2-7  &  {\green{ 1.0-6}}  & 03 & 07 & 15  \\ \hline 
 fap25 &  2118;0 &  2118 &   0;46 &  95 &  218 &  9.2-7  &  8.1-7  &  9.0-7  &  -3.1-8  &  5.1-8  &  {\green{ 1.1-6}}  & 2:44 & 11:02 & 12:26  \\ \hline 
 fap36 &  4110;0 &  4110 &  24;48 &  129 &  250 &  8.4-7  &  9.7-7  &  9.9-7  &  1.4-9  &  2.7-8  &  {\green{ 1.3-6}}  &   1:15:31 &   1:38:55 &   1:40:50  \\ \hline 
 bur26a &  1051;0 &   676 &  151;27 &  3110 &  43852 &  9.9-7  &  9.9-7  &  9.9-7  &    {\blue{ -5.7-6}}  &    {\blue{ -8.0-6}}  &    {\blue{ -5.6-6}}  & 6:14 & 17:12 &   1:33:31  \\ \hline 
 bur26b &  1051;0 &   676 &  150;27 &  3083 &  43600 &  9.6-7  &  9.9-7  &  9.9-7  &    {\blue{ -5.7-6}}  &    {\blue{ -7.9-6}}  &    {\blue{ -5.5-6}}  & 6:05 & 17:05 &   1:38:59  \\ \hline 
 bur26c &  1051;0 &   676 &  146;27 &  3245 &  45432 &  9.3-7  &  9.9-7  &  9.9-7  &    {\blue{ -6.9-6}}  &    {\blue{ -9.2-6}}  &    {\blue{ -6.3-6}}  & 6:18 & 18:47 &   2:08:59  \\ \hline 
 bur26d &  1051;0 &   676 &  172;27 &  3188 &  44790 &  9.6-7  &  9.9-7  &  9.9-7  &    {\blue{ -6.1-6}}  &    {\blue{ -8.7-6}}  &    {\blue{ -6.0-6}}  & 7:12 & 18:00 &   2:30:30  \\ \hline 
 bur26e &  1051;0 &   676 &  153;27 &  3178 &  44615 &  9.6-7  &  9.9-7  &  9.9-7  &    {\blue{ -6.2-6}}  &    {\blue{ -8.7-6}}  &    {\blue{ -6.0-6}}  & 8:01 & 18:23 &   2:49:49  \\ \hline 
 bur26f &  1051;0 &   676 &  127;27 &  3142 &  44272 &  9.9-7  &  9.9-7  &  9.9-7  &    {\blue{ -5.7-6}}  &    {\blue{ -8.6-6}}  &    {\blue{ -5.9-6}}  & 4:10 & 17:37 &   2:48:14  \\ \hline 
 bur26g &  1051;0 &   676 &  149;27 &  3012 &  42291 &  9.8-7  &  9.9-7  &  9.9-7  &    {\blue{ -6.3-6}}  &    {\blue{ -8.4-6}}  &    {\blue{ -5.9-6}}  & 7:18 & 17:23 &   2:06:01  \\ \hline 
 bur26h &  1051;0 &   676 &  126;27 &  3004 &  42050 &  9.7-7  &  9.9-7  &  9.9-7  &    {\blue{ -6.1-6}}  &    {\blue{ -8.3-6}}  &    {\blue{ -5.9-6}}  & 4:09 & 15:40 &   2:27:21  \\ \hline 
 chr20a &  628;0 &   400 &  242;40 &  5548 &  75000 &  9.4-7  &  9.9-7  &  {\green{ 1.3-6}}  &    {\red{ -1.7-5}}  &    {\red{ -2.7-5}}  &    {\red{ -2.1-5}}  & 2:44 & 10:09 &   1:39:07  \\ \hline 
 chr20b &  628;0 &   400 &  203;33 &  4727 &  74892 &  9.9-7  &  9.9-7  &  9.9-7  &    {\red{ -1.5-5}}  &    {\red{ -2.5-5}}  &    {\red{ -1.7-5}}  & 2:41 & 8:34 &   1:32:59  \\ \hline 
 chr20c &  628;0 &   400 &  223;48 &  6738 &  75000 &  9.7-7  &  9.9-7  &  {\green{ 1.6-6}}  &    {\red{ -2.1-5}}  &    {\red{ -3.4-5}}  &    {\red{ -3.0-5}}  & 3:25 & 12:27 & 59:58  \\ \hline 
 chr22a &  757;0 &   484 &  390;42 &  4871 &  75000 &  9.4-7  &  9.9-7  &  {\green{ 1.0-6}}  &    {\red{ -1.3-5}}  &    {\red{ -2.1-5}}  &    {\red{ -1.5-5}}  & 6:10 & 13:04 &   1:44:43  \\ \hline 
 chr22b &  757;0 &   484 &  224;40 &  4430 &  68291 &  9.9-7  &  9.9-7  &  9.9-7  &    {\red{ -1.4-5}}  &    {\red{ -2.0-5}}  &    {\red{ -1.4-5}}  & 3:31 & 11:51 &   1:44:34  \\ \hline 
 chr25a &  973;0 &   625 &  234;36 &  4750 &  75000 &  9.9-7  &  9.9-7  &  {\green{ 1.0-6}}  &    {\red{ -1.3-5}}  &    {\red{ -2.3-5}}  &    {\red{ -1.6-5}}  & 8:08 & 21:30 &   4:00:41  \\ \hline 
 esc32a &  1582;0 &  1024 &  155;48 &  4950 &  75000 &  9.7-7  &  9.7-7  &  {\green{ 1.2-6}}  &    {\red{ -1.2-5}}  &    {\red{ -1.9-5}}  &    {\red{ -1.3-5}}  & 8:41 & 52:28 &   4:35:48  \\ \hline 
 esc32b &  1582;0 &  1024 &  152;39 &  6966 &  75000 &  9.0-7  &  9.9-7  &  {\green{ 1.2-6}}  &    {\red{ -1.4-5}}  &    {\red{ -2.0-5}}  &    {\red{ -1.5-5}}  & 10:26 &   1:30:30 &   4:50:54  \\ \hline 
 esc32c &  1582;0 &  1024 &  177;77 &  10626 &  75000 &  9.7-7  &  9.9-7  &  {\green{ 3.0-6}}  &    {\red{ -3.4-5}}  &    {\red{ -4.6-5}}  &    {\red{ -5.4-5}}  & 13:49 &   2:48:01 &   4:58:03  \\ \hline 
 esc32d &  1582;0 &  1024 &  184;71 &  12520 &  75000 &  9.8-7  &  9.9-7  &  {\green{ 2.7-6}}  &    {\red{ -3.1-5}}  &    {\red{ -4.1-5}}  &    {\red{ -4.7-5}}  & 28:56 &   3:05:06 &   5:16:11  \\ \hline 
 esc32e &  1582;0 &  1024 &  272;38 &  3245 &  45376 &  9.2-7  &  9.9-7  &  9.9-7  &  -2.7-7  &    {\blue{ -5.5-6}}  &  {\green{ -3.8-6}}  & 18:12 & 23:50 &   3:32:35  \\ \hline 
 esc32f &  1582;0 &  1024 &  272;38 &  3245 &  45376 &  9.2-7  &  9.9-7  &  9.9-7  &  -2.7-7  &    {\blue{ -5.5-6}}  &  {\green{ -3.8-6}}  & 18:51 & 25:13 &   3:51:30  \\ \hline 
 esc32g &  1582;0 &  1024 &  113;46 &  4071 &  63055 &  9.5-7  &  9.9-7  &  9.9-7  &  {\green{ -1.2-6}}  &    {\blue{ -9.5-6}}  &    {\blue{ -6.7-6}}  & 13:00 & 31:09 &   5:41:29  \\ \hline 
 esc32h &  1582;0 &  1024 &  145;66 &  9791 &  75000 &  9.5-7  &  9.9-7  &  {\green{ 2.7-6}}  &    {\red{ -3.0-5}}  &    {\red{ -4.4-5}}  &    {\red{ -5.0-5}}  & 20:40 &   1:41:10 &   10:49:54  \\ \hline 
 had20 &  628;0 &   400 &  153;41 &  5269 &  75000 &  9.6-7  &  9.9-7  &  {\green{ 1.2-6}}  &    {\blue{ -9.8-6}}  &    {\red{ -1.4-5}}  &    {\red{ -1.0-5}}  & 1:38 & 8:02 & 41:08  \\ \hline 
 kra30a &  1393;0 &   900 &  179;42 &  6177 &  75000 &  9.9-7  &  9.9-7  &  {\green{ 1.5-6}}  &  {\green{ -3.8-6}}  &    {\red{ -3.6-5}}  &    {\red{ -2.9-5}}  & 5:55 & 57:39 &   3:51:14  \\ \hline 
 kra30b &  1393;0 &   900 &  240;42 &  6216 &  75000 &  9.9-7  &  9.9-7  &  {\green{ 1.5-6}}  &    {\red{ -2.6-5}}  &    {\red{ -3.7-5}}  &    {\red{ -3.0-5}}  & 19:19 &   1:11:10 &   3:58:05  \\ \hline 
 kra32 &  1582;0 &  1024 &  336;44 &  6190 &  75000 &  9.6-7  &  9.9-7  &  {\green{ 1.4-6}}  &    {\red{ -1.9-5}}  &    {\red{ -2.9-5}}  &    {\red{ -2.4-5}}  & 40:34 &   1:47:07 &   6:24:11  \\ \hline 
 lipa20a &  628;0 &   400 &  125;40 &  2273 &  29318 &  9.3-7  &  9.9-7  &  9.9-7  &  {\green{ -1.6-6}}  &  {\green{ -3.6-6}}  &  {\green{ -2.4-6}}  & 1:04 & 3:43 & 20:31  \\ \hline 
 lipa20b &  628;0 &   400 &  164;38 &  4268 &  66682 &  9.4-7  &  9.9-7  &  9.9-7  &    {\blue{ -5.9-6}}  &    {\red{ -1.1-5}}  &    {\blue{ -7.2-6}}  & 1:24 & 3:30 & 47:49  \\ \hline 
 lipa30a &  1393;0 &   900 &  193;37 &  1980 &  22720 &  9.8-7  &  9.9-7  &  9.9-7  &  -8.3-7  &  {\green{ -2.3-6}}  &  {\green{ -1.5-6}}  & 9:28 & 10:01 &   1:30:56  \\ \hline 
 lipa30b &  1393;0 &   900 &  193;29 &  3441 &  50659 &  9.4-7  &  9.9-7  &  9.9-7  &  {\green{ -4.4-6}}  &    {\blue{ -7.1-6}}  &  {\green{ -4.8-6}}  & 11:30 & 17:46 &   3:52:30  \\ \hline 
 lipa40a &  2458;0 &  1600 &  83; 9 &  1646 &  18212 &  9.7-7  &  9.9-7  &  9.9-7  &  -6.2-7  &  {\green{ -1.8-6}}  &  {\green{ -1.2-6}}  & 21:13 & 37:53 &   8:40:03  \\ \hline 
 lipa40b &  2458;0 &  1600 &  274;33 &  3460 &  51010 &  9.6-7  &  9.9-7  &  9.9-7  &  {\green{ -3.7-6}}  &    {\blue{ -7.2-6}}  &  {\green{ -4.8-6}}  &   1:11:40 &   2:20:07 &   14:32:56  \\ \hline 
 nug20 &  628;0 &   400 &  159;43 &  6400 &  75000 &  9.8-7  &  9.9-7  &  {\green{ 1.5-6}}  &    {\red{ -1.9-5}}  &    {\red{ -2.7-5}}  &    {\red{ -2.3-5}}  & 1:49 & 9:28 &   1:35:21  \\ \hline 
 nug21 &  691;0 &   441 &  204;45 &  6713 &  75000 &  9.9-7  &  9.9-7  &  {\green{ 1.6-6}}  &    {\red{ -2.4-5}}  &    {\red{ -3.0-5}}  &    {\red{ -2.6-5}}  & 5:24 & 13:12 &   1:58:15  \\ \hline 
 nug22 &  757;0 &   484 &  409;48 &  7196 &  75000 &  9.7-7  &  9.9-7  &  {\green{ 1.8-6}}  &    {\red{ -2.5-5}}  &    {\red{ -3.6-5}}  &    {\red{ -3.2-5}}  & 8:03 & 17:04 &   2:19:02  \\ \hline 
 nug24 &  898;0 &   576 &  181;43 &  6400 &  75000 &  9.6-7  &  9.9-7  &  {\green{ 1.5-6}}  &    {\red{ -1.9-5}}  &    {\red{ -2.7-5}}  &    {\red{ -2.2-5}}  & 6:25 & 23:30 &   1:45:44  \\ \hline 
 nug25 &  973;0 &   625 &  240;42 &  6250 &  75000 &  9.8-7  &  9.9-7  &  {\green{ 1.5-6}}  &    {\red{ -1.8-5}}  &    {\red{ -2.5-5}}  &    {\red{ -2.0-5}}  & 3:39 & 30:34 &   2:06:53  \\ \hline 
 nug27 &  1132;0 &   729 &  287;44 &  5217 &  75000 &  9.9-7  &  9.9-7  &  {\green{ 1.1-6}}  &    {\red{ -1.9-5}}  &    {\red{ -2.7-5}}  &    {\red{ -2.0-5}}  & 5:42 & 36:24 &   2:57:17  \\ \hline 
 nug28 &  1216;0 &   784 &  277;44 &  5048 &  75000 &  9.8-7  &  9.9-7  &  {\green{ 1.1-6}}  &    {\red{ -1.8-5}}  &    {\red{ -2.5-5}}  &    {\red{ -1.8-5}}  & 7:23 & 38:44 &   3:42:27  \\ \hline 
 nug30 &  1393;0 &   900 &  210;41 &  6156 &  75000 &  9.6-7  &  9.9-7  &  {\green{ 1.4-6}}  &    {\red{ -1.6-5}}  &    {\red{ -2.2-5}}  &    {\red{ -1.8-5}}  & 9:58 &   1:11:51 &   7:20:31  \\ \hline 
 rou20 &  628;0 &   400 &  156;36 &  5308 &  75000 &  9.9-7  &  9.9-7  &  {\green{ 1.2-6}}  &    {\red{ -1.3-5}}  &    {\red{ -1.6-5}}  &    {\red{ -1.2-5}}  & 1:21 & 8:50 &   1:39:12  \\ \hline 
 scr20 &  628;0 &   400 &  198;39 &  5555 &  75000 &  9.8-7  &  9.9-7  &  {\green{ 1.3-6}}  &    {\red{ -1.7-5}}  &    {\red{ -2.9-5}}  &    {\red{ -2.2-5}}  & 1:49 & 5:03 &   1:32:58  \\ \hline 
 ste36a &  1996;0 &  1296 &  211;37 &  5494 &  75000 &  9.8-7  &  9.9-7  &  {\green{ 1.2-6}}  &    {\red{ -1.7-5}}  &    {\red{ -2.5-5}}  &    {\red{ -1.9-5}}  & 40:38 &   1:11:46 &   8:28:01  \\ \hline 
 ste36b &  1996;0 &  1296 &  231;48 &  6918 &  75000 &  9.6-7  &  9.9-7  &  {\green{ 1.7-6}}  &    {\red{ -2.9-5}}  &    {\red{ -3.9-5}}  &    {\red{ -3.5-5}}  &   1:05:38 &   1:57:46 &   9:39:13  \\ \hline 
 ste36c &  1996;0 &  1296 &  225;39 &  5779 &  75000 &  9.9-7  &  9.9-7  &  {\green{ 1.3-6}}  &    {\red{ -1.8-5}}  &    {\red{ -2.7-5}}  &    {\red{ -2.1-5}}  & 50:54 &   2:46:30 &   14:10:06  \\ \hline 
 tai20a &  628;0 &   400 &  164;31 &  5072 &  75000 &  9.8-7  &  9.9-7  &  {\green{ 1.1-6}}  &    {\red{ -1.1-5}}  &    {\red{ -1.5-5}}  &    {\red{ -1.1-5}}  & 1:23 & 8:37 &   1:35:41  \\ \hline 
 tai20b &  628;0 &   400 &  180;45 &  7271 &  75000 &  9.9-7  &  9.9-7  &  {\green{ 1.8-6}}  &    {\red{ -2.8-5}}  &    {\red{ -3.9-5}}  &    {\red{ -3.6-5}}  & 3:26 & 5:52 &   1:27:02  \\ \hline 
 tai25a &  973;0 &   625 &  346;31 &  6297 &  75000 &  9.9-7  &  9.9-7  &  {\green{ 1.5-6}}  &  {\green{ -2.4-6}}  &    {\red{ -2.8-5}}  &    {\red{ -2.4-5}}  & 7:51 & 13:35 &   1:21:54  \\ \hline 
 tai25b &  973;0 &   625 &  141;39 &  6176 &  75000 &  9.9-7  &  9.9-7  &  {\green{ 1.4-6}}  &    {\red{ -1.8-5}}  &    {\red{ -2.6-5}}  &    {\red{ -2.2-5}}  & 8:02 & 13:42 &   1:23:59  \\ \hline 
 tai30a &  1393;0 &   900 &  258;31 &  4212 &  65722 &  9.9-7  &  9.9-7  &  9.9-7  &    {\blue{ -7.0-6}}  &    {\red{ -1.1-5}}  &    {\blue{ -7.5-6}}  & 17:39 & 19:37 &   2:43:13  \\ \hline 
 tai30b &  1393;0 &   900 &  106;35 &  4489 &  71044 &  9.9-7  &  9.9-7  &  9.9-7  &    {\blue{ -9.5-6}}  &    {\red{ -2.1-5}}  &    {\red{ -1.5-5}}  & 5:07 & 21:28 &   2:56:30  \\ \hline 
 tai35a &  1888;0 &  1225 &  246;30 &  3953 &  59536 &  9.7-7  &  9.9-7  &  9.9-7  &    {\blue{ -5.9-6}}  &    {\blue{ -8.9-6}}  &    {\blue{ -6.0-6}}  & 13:08 & 40:30 &   5:06:41  \\ \hline 
 tai35b &  1888;0 &  1225 &  128;43 &  6624 &  75000 &  9.9-7  &  9.9-7  &  {\green{ 1.6-6}}  &    {\red{ -1.7-5}}  &    {\red{ -2.4-5}}  &    {\red{ -2.1-5}}  & 19:14 &   1:09:33 &   6:24:39  \\ \hline 
 tai40a &  2458;0 &  1600 &  325;31 &  3973 &  60071 &  9.8-7  &  9.9-7  &  9.9-7  &  {\green{ -5.0-6}}  &    {\blue{ -8.6-6}}  &    {\blue{ -5.7-6}}  & 32:43 &   1:24:16 &   9:43:53  \\ \hline 
 tai40b &  2458;0 &  1600 &  126;33 &  4316 &  66912 &  9.2-7  &  9.9-7  &  9.9-7  &    {\red{ -1.1-5}}  &    {\red{ -1.5-5}}  &    {\red{ -1.0-5}}  & 29:47 &   2:43:34 &   12:22:13  \\ \hline 
 tho30 &  1393;0 &   900 &  159;43 &  6624 &  75000 &  9.9-7  &  9.9-7  &  {\green{ 1.6-6}}  &    {\red{ -1.7-5}}  &    {\red{ -2.3-5}}  &    {\red{ -2.0-5}}  & 12:49 &   1:16:52 &   3:51:02  \\ \hline 
 tho40 &  2458;0 &  1600 &  194;43 &  6242 &  75000 &  9.8-7  &  9.9-7  &  {\green{ 1.5-6}}  &    {\red{ -1.6-5}}  &    {\red{ -2.1-5}}  &    {\red{ -1.8-5}}  &   1:41:24 &   3:51:16 &   24:09:34  \\ \hline 
 be250.1 &  251;0 &   251 &   0;707 &  2448 &  2547 &  9.9-7  &  9.9-7  &  9.1-7  &  4.5-7  &  {\green{ -1.2-6}}  &  {\green{ -2.6-6}}  & 24 & 1:05 & 39  \\ \hline 
 be250.2 &  251;0 &   251 &   0;693 &  2571 &  2525 &  9.8-7  &  9.9-7  &  9.5-7  &  6.5-7  &    {\blue{ -5.1-6}}  &  {\green{ -2.7-6}}  & 15 & 57 & 41  \\ \hline 
 be250.3 &  251;0 &   251 &   0;613 &  2346 &  2525 &  9.9-7  &  9.9-7  &  9.7-7  &  -4.4-7  &  2.5-7  &  {\green{ -2.8-6}}  & 22 & 1:11 & 41  \\ \hline 
 be250.4 &  251;0 &   251 &   0;694 &  2482 &  2547 &  9.9-7  &  9.9-7  &  8.5-7  &  2.7-7  &  {\green{ -4.5-6}}  &  {\green{ -2.7-6}}  & 23 & 1:14 & 41  \\ \hline 
 be250.5 &  251;0 &   251 &   0;598 &  878 &  2598 &  9.9-7  &  9.9-7  &  9.4-7  &  -2.8-7  &  {\green{ -2.0-6}}  &  {\green{ 1.6-6}}  & 26 & 28 & 42  \\ \hline 
 be250.6 &  251;0 &   251 &   0;599 &  893 &  2547 &  9.9-7  &  9.9-7  &  8.9-7  &  -5.7-7  &    {\blue{ -5.7-6}}  &  {\green{ -2.5-6}}  & 14 & 24 & 40  \\ \hline 
 be250.7 &  251;0 &   251 &   0;621 &  2959 &  2525 &  9.9-7  &  9.9-7  &  9.6-7  &  -4.8-7  &  {\green{ 3.4-6}}  &  {\green{ -2.9-6}}  & 21 & 1:26 & 41  \\ \hline 
 be250.8 &  251;0 &   251 &   0;649 &  995 &  2547 &  9.9-7  &  9.9-7  &  8.5-7  &  -8.1-8  &  8.5-7  &  {\green{ -2.7-6}}  & 26 & 29 & 40  \\ \hline 
 be250.9 &  251;0 &   251 &   0;743 &  3326 &  2546 &  9.9-7  &  9.9-7  &  9.9-7  &  4.5-7  &  -9.4-7  &  -8.3-8  & 24 & 1:30 & 41  \\ \hline 
 be250.10 &  251;0 &   251 &   0;637 &  2716 &  2547 &  9.9-7  &  9.9-7  &  8.9-7  &  -1.4-7  &  {\green{ 4.6-6}}  &  {\green{ -2.7-6}}  & 15 & 1:10 & 40  \\ \hline 
 bqp500-1 &  501;0 &   501 &   0;676 &  2378 &  3365 &  9.9-7  &  9.7-7  &  9.3-7  &  {\green{ 1.2-6}}  &    {\blue{ 9.0-6}}  &  9.0-7  & 2:59 & 8:19 & 7:23  \\ \hline 
 bqp500-2 &  501;0 &   501 &   0;741 &  2562 &  4448 &  9.2-7  &  9.8-7  &  9.9-7  &  {\green{ 1.0-6}}  &    {\blue{ 9.9-6}}  &    {\blue{ -5.3-6}}  & 3:15 & 9:08 & 9:54  \\ \hline 
 bqp500-3 &  501;0 &   501 &   0;674 &  2550 &  4430 &  9.9-7  &  9.9-7  &  9.9-7  &  {\green{ 1.4-6}}  &    {\red{ 1.0-5}}  &    {\blue{ -5.4-6}}  & 3:01 & 9:01 & 9:53  \\ \hline 
 bqp500-4 &  501;0 &   501 &   0;648 &  2574 &  4503 &  9.9-7  &  9.9-7  &  9.9-7  &  {\green{ 1.2-6}}  &    {\red{ 1.1-5}}  &    {\blue{ -5.5-6}}  & 2:53 & 9:00 & 10:09  \\ \hline 
 bqp500-5 &  501;0 &   501 &   0;760 &  2507 &  4000 &  9.9-7  &  9.9-7  &  9.9-7  &  8.7-7  &    {\blue{ 9.7-6}}  &    {\blue{ -5.2-6}}  & 3:35 & 8:48 & 8:41  \\ \hline 
 bqp500-6 &  501;0 &   501 &   0;733 &  2450 &  3472 &  9.9-7  &  9.9-7  &  9.6-7  &  8.4-7  &    {\blue{ 9.3-6}}  &  {\green{ -3.3-6}}  & 3:17 & 8:26 & 7:19  \\ \hline 
 bqp500-7 &  501;0 &   501 &   0;709 &  2483 &  3621 &  9.9-7  &  9.9-7  &  9.9-7  &  9.9-7  &    {\blue{ 9.8-6}}  &    {\blue{ -5.0-6}}  & 2:56 & 8:05 & 7:34  \\ \hline 
 bqp500-8 &  501;0 &   501 &   0;713 &  2449 &  3472 &  9.8-7  &  9.9-7  &  9.3-7  &  {\green{ 1.4-6}}  &    {\blue{ 9.5-6}}  &  {\green{ -3.6-6}}  & 2:28 & 7:43 & 7:19  \\ \hline 
 bqp500-9 &  501;0 &   501 &   0;674 &  2452 &  3431 &  9.9-7  &  9.8-7  &  9.8-7  &  7.7-7  &    {\blue{ 9.6-6}}  &  {\green{ -3.3-6}}  & 3:02 & 8:32 & 7:12  \\ \hline 
 bqp500-10 &  501;0 &   501 &   0;692 &  2570 &  4477 &  9.9-7  &  9.9-7  &  9.9-7  &  {\green{ 1.1-6}}  &    {\red{ 1.0-5}}  &    {\blue{ -5.4-6}}  & 2:47 & 8:22 & 8:16  \\ \hline 

 gka1e &  201;0 &   201 &   0;675 &  958 &  2999 &  9.9-7  &  9.9-7  &  9.7-7  &  7.2-7  &  {\green{ 4.8-6}}  &  {\green{ 2.8-6}}  & 18 & 19 & 31  \\ \hline 
 gka2e &  201;0 &   201 &   0;719 &  1155 &  2632 &  9.9-7  &  9.9-7  &  9.9-7  &  9.1-7  &    {\blue{ 5.9-6}}  &  9.2-7  & 19 & 20 & 27  \\ \hline 
 gka3e &  201;0 &   201 &   0;753 &  2995 &  3730 &  9.9-7  &  9.9-7  &  9.6-7  &  5.0-7  &  {\green{ -2.4-6}}  &  {\green{ 2.8-6}}  & 22 & 1:01 & 37  \\ \hline 
 gka4e &  201;0 &   201 &   0;735 &  2445 &  2036 &  9.9-7  &  9.9-7  &  9.9-7  &  -2.5-7  &  {\green{ -1.7-6}}  &  -8.0-7  & 21 & 47 & 21  \\ \hline 
 gka5e &  201;0 &   201 &   0;712 &  2931 &  2598 &  9.9-7  &  9.9-7  &  9.8-7  &  4.6-8  &  {\green{ -3.5-6}}  &  {\green{ -3.4-6}}  & 19 & 57 & 27  \\ \hline 
 gka1f &  501;0 &   501 &   0;630 &  2233 &  2731 &  9.8-7  &  9.9-7  &  9.8-7  &  4.4-7  &    {\blue{ -7.4-6}}  &  {\green{ -4.2-6}}  & 2:53 & 7:55 & 5:54  \\ \hline 
 gka2f &  501;0 &   501 &   0;684 &  2298 &  3338 &  9.7-7  &  9.8-7  &  9.9-7  &  8.4-7  &    {\blue{ 9.1-6}}  &  {\green{ 3.2-6}}  & 3:05 & 8:11 & 7:11  \\ \hline 
 gka3f &  501;0 &   501 &   0;801 &  2579 &  4430 &  9.9-7  &  9.9-7  &  9.9-7  &  9.3-7  &    {\red{ 1.0-5}}  &    {\blue{ -5.4-6}}  & 3:29 & 8:40 & 9:39  \\ \hline 
 gka4f &  501;0 &   501 &   0;751 &  2798 &  5158 &  9.9-7  &  9.9-7  &  9.8-7  &  8.3-7  &    {\red{ 1.1-5}}  &  {\green{ 4.2-6}}  & 2:36 & 9:15 & 11:20  \\ \hline 
 gka5f &  501;0 &   501 &   0;749 &  2120 &  5158 &  9.7-7  &  9.9-7  &  9.9-7  &  {\green{ 1.0-6}}  &    {\blue{ -9.6-6}}  &  {\green{ 4.2-6}}  & 3:28 & 7:20 & 9:43  \\ \hline 
 soybean-large.2 &  308;0 &   307 &   0;1745 &  6307 &  13795 &  9.8-7  &  9.8-7  &  9.7-7  &  -7.2-7  &  -9.4-10  &  -4.9-8  & 1:17 & 3:07 & 5:03  \\ \hline 
 soybean-large.3 &  308;0 &   307 &   0;1764 &  4676 &  14096 &  9.1-7  &  9.9-7  &  7.3-7  &  -1.8-7  &  -3.8-8  &  -1.2-8  & 45 & 2:17 & 5:22  \\ \hline 
 soybean-large.4 &  308;0 &   307 &   0;1262 &  4417 &  6965 &  8.5-7  &  9.9-7  &  8.9-7  &  1.2-7  &  7.2-9  &  -1.9-8  & 34 & 2:16 & 2:31  \\ \hline 
 soybean-large.5 &  308;0 &   307 &   0;1017 &  4492 &  6219 &  9.4-7  &  9.9-7  &  9.8-7  &  -2.9-8  &  2.4-9  &  -2.2-8  & 43 & 2:11 & 2:19  \\ \hline 
 soybean-large.6 &  308;0 &   307 &   0;601 &  2623 &  7468 &  9.8-7  &  9.9-7  &  9.8-7  &  -8.3-9  &  -1.6-9  &  -3.7-11  & 26 & 1:19 & 2:45  \\ \hline 
 soybean-large.7 &  308;0 &   307 &   0;740 &  2103 &  4065 &  9.8-7  &  9.9-7  &  8.7-7  &  5.8-9  &  -7.0-10  &  1.0-9  & 32 & 1:03 & 1:26  \\ \hline 
 soybean-large.8 &  308;0 &   307 &   0;706 &  2501 &  4351 &  9.9-7  &  9.9-7  &  9.9-7  &  -1.3-8  &  2.9-9  &  -8.4-9  & 34 & 1:18 & 1:34  \\ \hline 
 soybean-large.9 &  308;0 &   307 &   0;728 &  2289 &  4560 &  9.4-7  &  9.9-7  &  7.8-7  &  3.1-8  &  -2.4-9  &  2.7-8  & 27 & 1:20 & 1:43  \\ \hline 
 soybean-large.10 &  308;0 &   307 &   0;525 &  2227 &  4560 &  9.5-7  &  9.9-7  &  9.3-7  &  5.5-8  &  -2.5-9  &  3.0-8  & 15 & 1:15 & 1:39  \\ \hline 
 soybean-large.11 &  308;0 &   307 &   0;484 &  3806 &  2261 &  9.7-7  &  9.9-7  &  8.3-7  &  1.7-8  &  5.4-9  &  -1.3-8  & 17 & 1:57 & 46  \\ \hline 

 spambase-large.2 &  1501;0 &  1500 &   0;62 &  217 &  263 &  5.4-7  &  8.2-7  &  2.8-8  &  -7.0-8  &  -4.0-10  &  4.0-10  & 4:02 & 8:40 & 6:47  \\ \hline 
 spambase-large.3 &  1501;0 &  1500 &   0;62 &  217 &  263 &  5.4-7  &  8.2-7  &  2.8-8  &  -7.0-8  &  -3.7-10  &  4.2-10  & 3:07 & 7:36 & 6:47  \\ \hline 
 spambase-large.4 &  1501;0 &  1500 &   0;62 &  217 &  263 &  5.4-7  &  8.2-7  &  2.8-8  &  -7.0-8  &  -3.5-10  &  4.4-10  & 3:47 & 8:39 & 6:47  \\ \hline 
 spambase-large.5 &  1501;0 &  1500 &   0;62 &  217 &  263 &  5.4-7  &  8.2-7  &  2.8-8  &  -7.0-8  &  -3.2-10  &  4.5-10  & 3:53 & 9:30 & 6:38  \\ \hline 
 spambase-large.6 &  1501;0 &  1500 &   0;62 &  217 &  263 &  5.4-7  &  8.2-7  &  2.9-8  &  -7.0-8  &  -3.0-10  &  4.7-10  & 3:41 & 9:34 & 6:22  \\ \hline 
 spambase-large.7 &  1501;0 &  1500 &   0;62 &  217 &  263 &  5.4-7  &  8.2-7  &  2.9-8  &  -7.0-8  &  -2.7-10  &  4.9-10  & 3:47 & 10:21 & 6:36  \\ \hline 
 spambase-large.8 &  1501;0 &  1500 &   0;62 &  217 &  263 &  5.4-7  &  8.2-7  &  2.9-8  &  -7.0-8  &  -2.5-10  &  5.0-10  & 3:56 & 9:55 & 6:40  \\ \hline 
 spambase-large.9 &  1501;0 &  1500 &   0;62 &  217 &  263 &  5.4-7  &  8.2-7  &  2.9-8  &  -7.0-8  &  -2.2-10  &  5.2-10  & 3:53 & 9:14 & 6:19  \\ \hline 
 spambase-large.10 &  1501;0 &  1500 &   0;62 &  217 &  263 &  5.4-7  &  8.2-7  &  2.9-8  &  -7.0-8  &  -2.0-10  &  5.4-10  & 4:12 & 8:59 & 5:36  \\ \hline 
 spambase-large.11 &  1501;0 &  1500 &   0;62 &  217 &  263 &  5.4-7  &  8.2-7  &  2.9-8  &  -7.0-8  &  -1.7-10  &  5.5-10  & 4:08 & 8:04 & 7:12  \\ \hline 

 abalone-large.2 &  1001;0 &  1000 &   0;869 &  10956 &  80482 &  9.9-7  &  9.8-7  &  9.9-7  &  1.1-8  &  3.4-9  &  -1.3-8  & 11:45 &   2:26:39 &   8:46:18  \\ \hline 
 abalone-large.3 &  1001;0 &  1000 &   0;695 &  5033 &  41764 &  9.9-7  &  9.9-7  &  9.9-7  &  -2.2-8  &  -4.2-9  &  -4.0-9  & 10:44 &   1:10:46 &   4:53:35  \\ \hline 
 abalone-large.4 &  1001;0 &  1000 &   0;716 &  3089 &  18374 &  8.2-7  &  9.9-7  &  6.9-7  &  5.9-8  &  -1.3-9  &  -2.6-10  & 10:49 & 41:57 &   2:24:14  \\ \hline 
 abalone-large.5 &  1001;0 &  1000 &   0;590 &  2872 &  10594 &  9.6-7  &  9.9-7  &  6.0-7  &  4.6-8  &  -6.5-10  &  -1.9-9  & 9:28 & 39:26 &   1:33:36  \\ \hline 
 abalone-large.6 &  1001;0 &  1000 &   0;561 &  3104 &  7697 &  9.7-7  &  9.9-7  &  9.3-7  &  1.0-8  &  -1.9-9  &  -1.5-9  & 9:01 & 43:42 &   1:07:38  \\ \hline 
 abalone-large.7 &  1001;0 &  1000 &   0;522 &  2949 &  7470 &  9.7-7  &  9.9-7  &  9.9-7  &  -5.5-9  &  8.8-10  &  2.6-9  & 8:43 & 40:16 &   1:07:55  \\ \hline 
 abalone-large.8 &  1001;0 &  1000 &   0;490 &  2920 &  5533 &  9.7-7  &  9.9-7  &  9.5-7  &  -1.9-8  &  1.8-9  &  4.0-9  & 7:53 & 36:26 & 51:42  \\ \hline 
 abalone-large.9 &  1001;0 &  1000 &   0;507 &  3117 &  3509 &  9.9-7  &  9.9-7  &  8.7-7  &  -1.9-10  &  -2.5-9  &  6.4-9  & 8:47 & 39:25 & 31:59  \\ \hline 
 abalone-large.10 &  1001;0 &  1000 &   0;503 &  3078 &  3202 &  9.9-7  &  9.9-7  &  8.6-7  &  -6.6-9  &  3.2-9  &  5.7-10  & 9:21 & 35:18 & 29:41  \\ \hline 
 abalone-large.11 &  1001;0 &  1000 &   0;531 &  2578 &  3206 &  9.5-7  &  9.9-7  &  9.3-7  &  3.0-8  &  2.0-9  &  6.0-10  & 8:26 & 32:41 & 28:18  \\ \hline 

 segment-large.2 &  1001;0 &  1000 &   0;52 &  121 &  210 &  8.8-7  &  9.7-7  &  7.2-7  &  -2.0-7  &  6.6-9  &  1.6-8  & 1:09 & 1:38 & 1:57  \\ \hline 
 segment-large.3 &  1001;0 &  1000 &   0;52 &  121 &  210 &  8.8-7  &  9.7-7  &  7.3-7  &  -2.0-7  &  6.8-9  &  1.6-8  & 1:13 & 1:44 & 1:55  \\ \hline 
 segment-large.4 &  1001;0 &  1000 &   0;52 &  121 &  210 &  8.8-7  &  9.7-7  &  7.3-7  &  -2.0-7  &  7.1-9  &  1.7-8  & 1:04 & 1:35 & 1:56  \\ \hline 
 segment-large.5 &  1001;0 &  1000 &   0;52 &  121 &  210 &  8.8-7  &  9.7-7  &  7.3-7  &  -2.0-7  &  7.4-9  &  1.7-8  & 1:08 & 1:31 & 1:54  \\ \hline 
 segment-large.6 &  1001;0 &  1000 &   0;52 &  121 &  210 &  8.8-7  &  9.7-7  &  7.4-7  &  -2.0-7  &  7.6-9  &  1.7-8  & 1:04 & 1:34 & 1:55  \\ \hline 
 segment-large.7 &  1001;0 &  1000 &   0;52 &  121 &  210 &  8.8-7  &  9.7-7  &  7.4-7  &  -2.0-7  &  7.9-9  &  1.7-8  & 1:05 & 1:32 & 1:46  \\ \hline 
 segment-large.8 &  1001;0 &  1000 &   0;52 &  121 &  210 &  8.8-7  &  9.7-7  &  7.5-7  &  -2.0-7  &  8.2-9  &  1.7-8  & 55 & 1:33 & 1:46  \\ \hline 
 segment-large.9 &  1001;0 &  1000 &   0;52 &  121 &  210 &  8.8-7  &  9.7-7  &  7.5-7  &  -2.0-7  &  8.4-9  &  1.7-8  & 48 & 1:30 & 1:42  \\ \hline 
 segment-large.10 &  1001;0 &  1000 &   0;52 &  121 &  210 &  8.8-7  &  9.7-7  &  7.6-7  &  -2.0-7  &  8.7-9  &  1.8-8  & 46 & 1:23 & 1:37  \\ \hline 
 segment-large.11 &  1001;0 &  1000 &   0;52 &  121 &  210 &  8.8-7  &  9.7-7  &  7.7-7  &  -2.0-7  &  8.9-9  &  1.8-8  & 43 & 1:40 & 1:27  \\ \hline 
 housing.2 &  507;0 &   506 &   0;15 &  47 &  54 &  1.8-7  &  8.2-7  &  4.7-7  &  3.4-9  &  -3.5-8  &  4.1-9  & 03 & 07 & 04  \\ \hline 
 housing.3 &  507;0 &   506 &   0;15 &  47 &  54 &  1.8-7  &  8.2-7  &  4.7-7  &  3.4-9  &  -3.5-8  &  4.1-9  & 03 & 05 & 04  \\ \hline 
 housing.4 &  507;0 &   506 &   0;15 &  47 &  54 &  1.7-7  &  8.1-7  &  4.7-7  &  3.4-9  &  -3.5-8  &  4.2-9  & 03 & 06 & 04  \\ \hline 
 housing.5 &  507;0 &   506 &   0;15 &  47 &  54 &  1.7-7  &  8.1-7  &  4.7-7  &  3.4-9  &  -3.5-8  &  4.2-9  & 03 & 08 & 04  \\ \hline 
 housing.6 &  507;0 &   506 &   0;15 &  47 &  54 &  1.7-7  &  8.1-7  &  4.7-7  &  3.4-9  &  -3.5-8  &  4.2-9  & 03 & 08 & 04  \\ \hline 
 housing.7 &  507;0 &   506 &   0;15 &  47 &  54 &  1.7-7  &  8.1-7  &  4.7-7  &  3.4-9  &  -3.5-8  &  4.2-9  & 02 & 08 & 04  \\ \hline 
 housing.8 &  507;0 &   506 &   0;15 &  47 &  54 &  1.6-7  &  8.1-7  &  4.7-7  &  3.4-9  &  -3.5-8  &  4.2-9  & 03 & 07 & 06  \\ \hline 
 housing.9 &  507;0 &   506 &   0;15 &  47 &  54 &  1.6-7  &  8.1-7  &  4.7-7  &  3.4-9  &  -3.5-8  &  4.2-9  & 04 & 08 & 06  \\ \hline 
 housing.10 &  507;0 &   506 &   0;15 &  47 &  54 &  1.6-7  &  8.1-7  &  4.7-7  &  3.4-9  &  -3.5-8  &  4.2-9  & 04 & 08 & 06  \\ \hline 
 housing.11 &  507;0 &   506 &   0;15 &  47 &  54 &  1.6-7  &  8.1-7  &  4.7-7  &  3.4-9  &  -3.5-8  &  4.2-9  & 03 & 08 & 07  \\ \hline 

 ex-be250.1 &  251;93375 &   251 &   0;719 &  9100 &  3005 &  9.9-7  &  9.9-7  &  9.9-7  &  -4.1-8  &  -6.1-9  &  -1.7-8  & 1:16 & 14:51 & 2:38  \\ \hline 
 ex-be250.2 &  251;93375 &   251 &   0;705 &  9409 &  2656 &  9.9-7  &  9.9-7  &  9.8-7  &  -7.8-8  &  1.4-9  &  -3.0-8  & 1:25 & 14:22 & 2:09  \\ \hline 
 ex-be250.3 &  251;93375 &   251 &   0;706 &  9195 &  2733 &  9.9-7  &  9.9-7  &  9.9-7  &  -3.4-8  &  -2.2-9  &  2.1-9  & 1:10 & 10:56 & 2:29  \\ \hline 
 ex-be250.4 &  251;93375 &   251 &   0;705 &  12298 &  2976 &  9.9-7  &  9.9-7  &  9.7-7  &  -5.0-8  &  -7.3-9  &  -8.8-9  & 1:12 & 17:06 & 2:29  \\ \hline 
 ex-be250.5 &  251;93375 &   251 &   0;705 &  10136 &  2894 &  9.8-7  &  9.9-7  &  9.9-7  &  -8.6-9  &  -1.4-8  &  -3.5-8  & 1:13 & 12:38 & 2:16  \\ \hline 
 ex-be250.6 &  251;93375 &   251 &   0;685 &  10050 &  2701 &  9.9-7  &  9.9-7  &  9.9-7  &  -5.7-8  &  4.0-9  &  -2.7-8  & 1:10 & 16:06 & 2:10  \\ \hline 
 ex-be250.7 &  251;93375 &   251 &   0;711 &  10887 &  2892 &  9.9-7  &  9.9-7  &  9.9-7  &  -2.9-8  &  -5.2-9  &  -1.2-8  & 1:12 & 16:42 & 2:19  \\ \hline 
 ex-be250.8 &  251;93375 &   251 &   0;712 &  9672 &  2733 &  9.9-7  &  9.9-7  &  9.8-7  &  -4.6-8  &  6.7-9  &  -2.8-8  & 1:13 & 12:50 & 2:17  \\ \hline 
 ex-be250.9 &  251;93375 &   251 &   0;711 &  6718 &  2645 &  9.9-7  &  9.9-7  &  9.8-7  &  -6.8-8  &  5.9-9  &  -5.1-9  & 1:18 & 10:21 & 2:05  \\ \hline 
 ex-be250.10 &  251;93375 &   251 &   0;717 &  11436 &  2843 &  9.9-7  &  9.9-7  &  9.7-7  &  -6.9-8  &  -1.9-9  &  -3.4-9  & 1:16 & 14:40 & 2:27  \\ \hline 

 ex-bqp500-1 &  501;374250 &   501 &   0;1439 &  15676 &  16436 &  9.9-7  &  9.9-7  &  8.8-7  &  -7.8-7  &  3.1-8  &  -4.1-7  & 14:34 &   2:00:25 &   1:42:18  \\ \hline 
 ex-bqp500-2 &  501;374250 &   501 &   0;1545 &  14458 &  17465 &  9.9-7  &  9.9-7  &  9.6-7  &  -8.2-7  &  4.4-8  &  -3.7-7  & 16:01 & 58:40 &   1:56:40  \\ \hline 
 ex-bqp500-3 &  501;374250 &   501 &   0;1795 &  15400 &  19460 &  9.9-7  &  9.9-7  &  9.5-7  &  -9.0-7  &  4.4-8  &  -2.5-7  & 19:43 &   1:12:55 &   1:51:14  \\ \hline 
 ex-bqp500-4 &  501;374250 &   501 &   0;1738 &  20092 &  19899 &  9.9-7  &  9.9-7  &  9.9-7  &  -9.1-7  &  4.8-8  &  -4.6-7  & 12:39 &   2:25:38 &   2:07:29  \\ \hline 
 ex-bqp500-5 &  501;374250 &   501 &   0;1777 &  15230 &  19517 &  9.9-7  &  9.9-7  &  9.3-7  &  -9.0-7  &  2.3-8  &  -2.2-7  & 17:38 &   1:01:55 &   2:10:26  \\ \hline 
 ex-bqp500-6 &  501;374250 &   501 &   0;1902 &  16149 &  17590 &  9.5-7  &  9.9-7  &  9.0-7  &  -7.0-7  &  5.9-8  &  -5.0-7  & 19:49 &   1:20:43 &   1:36:51  \\ \hline 
 ex-bqp500-7 &  501;374250 &   501 &   0;1835 &  15655 &  17481 &  9.9-7  &  9.9-7  &  9.9-7  &  -8.1-7  &  -6.2-8  &  -4.5-7  & 20:18 &   1:57:19 &   1:51:46  \\ \hline 
 ex-bqp500-8 &  501;374250 &   501 &   0;1798 &  16958 &  22473 &  9.9-7  &  9.9-7  &  9.7-7  &  -8.8-7  &  -8.2-9  &  -1.9-7  & 18:12 &   1:09:32 &   2:31:37  \\ \hline 
 ex-bqp500-9 &  501;374250 &   501 &   0;1977 &  15086 &  19588 &  9.9-7  &  9.9-7  &  9.9-7  &  -6.6-7  &  5.0-8  &  -2.6-7  & 20:19 &   1:14:18 &   2:04:43  \\ \hline 
 ex-bqp500-10 &  501;374250 &   501 &   0;1447 &  16325 &  19686 &  9.7-7  &  9.9-7  &  9.8-7  &  -7.6-7  &  -6.2-8  &  -2.9-7  & 15:30 &   2:01:24 &   2:06:46  \\ \hline 

ex-gka1e &  201;59700 &   201 &   0;750 &  7944 &  3909 &  9.9-7  &  9.9-7  &  9.9-7  &  -2.6-8  &  7.6-9  &  -4.1-9  & 35 & 7:18 & 1:38  \\ \hline 
 ex-gka2e &  201;59700 &   201 &   0;1651 &  25000 &  9885 &  9.9-7  &    {\red{ 1.0-5}}  &  9.6-7  &  5.9-7  &  1.0-7  &  -3.7-7  & 1:17 & 19:30 & 3:53  \\ \hline 
 ex-gka3e &  201;59700 &   201 &   0;1577 &  15270 &  13640 &  9.9-7  &  9.9-7  &  9.2-7  &  -3.2-7  &  -9.9-9  &  -3.3-7  & 1:07 & 14:41 & 5:02  \\ \hline 
 ex-gka4e &  201;59700 &   201 &   0;1460 &  25000 &  19976 &  9.9-7  &    {\red{ 1.7-5}}  &  9.9-7  &  -6.0-7  &  3.3-9  &  -5.2-7  & 1:20 & 22:57 & 7:34  \\ \hline 
 ex-gka5e &  201;59700 &   201 &   0;2038 &  25000 &  20715 &  9.9-7  &    {\red{ 1.9-5}}  &  9.9-7  &  -9.0-7  &  8.4-8  &  -6.1-7  & 1:16 & 22:12 & 7:42  \\ \hline 
 ex-gka1f &  501;374250 &   501 &   0;643 &  5599 &  2864 &  9.9-7  &  9.9-7  &  9.7-7  &  -8.4-8  &  -3.8-9  &  -1.4-9  & 7:15 & 50:41 & 18:28  \\ \hline 
 ex-gka2f &  501;374250 &   501 &   0;1049 &  8569 &  8939 &  9.9-7  &  9.9-7  &  9.8-7  &  -4.1-7  &  3.7-8  &  -1.8-7  & 6:18 &   1:08:19 & 55:07  \\ \hline 
 ex-gka3f &  501;374250 &   501 &   0;1798 &  21403 &  26462 &  9.9-7  &  9.9-7  &  9.2-7  &  -9.4-7  &  -5.6-8  &  -5.9-7  & 13:14 &   2:12:26 &   2:53:59  \\ \hline 
 ex-gka4f &  501;374250 &   501 &   0;2169 &  25000 &  37127 &  9.9-7  &  {\green{ 2.2-6}}  &  9.9-7  &  -7.2-7  &  -2.2-8  &  -6.9-7  & 21:49 &   2:27:06 &   3:11:07  \\ \hline 
 ex-gka5f &  501;374250 &   501 &   0;2259 &  25000 &  37456 &  9.8-7  &  {\green{ 1.6-6}}  &  9.9-7  &  {\green{ -1.2-6}}  &  -5.5-8  &  {\green{ -1.1-6}}  & 24:17 &   2:26:09 &   4:00:00  \\ \hline 
 		
		\end{longtable}
		
	\end{tiny}
\end{center}
\end{landscape}

\end{document}